\documentclass[11pt]{article}

\usepackage[dvipsnames]{xcolor}
\definecolor{SNSRed}{RGB}{183,42,58}
\definecolor{SNSBlue}{RGB}{0,114,142}
\definecolor{verdelink}{RGB}{6,168,21}

\usepackage{geometry}
\usepackage{calc}
\usepackage[colorlinks=true,hyperindex,linkcolor=SNSRed, pagebackref=false, citecolor=SNSBlue, urlcolor=verdelink, pdfpagelabels]{hyperref}

\usepackage{amsxtra}
\usepackage{amscd}
\usepackage{amsthm}
\usepackage{amsfonts}
\usepackage{amssymb}
\usepackage{eucal}
\usepackage{tikz-cd}
\usepackage[matrix,arrow,curve]{xy}
\usepackage{cleveref}
\usepackage[T1]{fontenc}
\usepackage[utf8]{inputenc}
\usepackage{mathtools}
\usepackage{float}

\theoremstyle{plain}
\newtheorem{Thm}{Theorem}[section]
\newtheorem{Cor}[Thm]{Corollary}
\newtheorem{Lem}[Thm]{Lemma}
\newtheorem{Conj}[Thm]{Conjecture}
\newtheorem{Prop}[Thm]{Proposition}

\theoremstyle{definition}
\newtheorem{Def}[Thm]{Definition}
\newtheorem{Ex}[Thm]{Example}

\theoremstyle{remark}
\newtheorem{Rem}[Thm]{Remark}

\newtheorem{conj}{Conjecture}

\numberwithin{equation}{section}

\renewcommand{\sec}[2]{\section{#2}\label{S:#1}}
\newcommand{\ssec}[2]{\subsection{#2}\label{SubS:#1}}
\newcommand{\sssec}[2]{\subsubsection{#2}\label{SubSubS:#1}}

 \newenvironment{thm}[1]%
    { \begin{Thm} \label{T:#1}}%
    { \end{Thm} }
\renewcommand{\th}[1]{\begin{thm}{#1} \sl }
\renewcommand{\eth}{\end{thm} }

\newenvironment{thms}[1]%
    { \begin{ThmS} \label{T:#1}}%
    { \end{ThmS} }
\newcommand{\ths}[1]{\begin{thms}{#1} \sl }
\newcommand{\eths}{\end{thms} }

\newenvironment{lemma}[1]%
    { \begin{Lem} \label{L:#1}}%
    { \end{Lem} }
\newcommand{\lem}[1]{\begin{lemma}{#1} \sl}
\newcommand{\elem}{\end{lemma}}

\newenvironment{lemmas}[1]%
    { \begin{LemS} \label{L:#1}}%
    { \end{LemS} }
\newcommand{\lems}[1]{\begin{lemmas}{#1} \sl}
\newcommand{\elems}{\end{lemmas}}

\newenvironment{propos}[1]%
    { \begin{Prop} \label{P:#1}}%
    { \end{Prop} }
\newcommand{\prop}[1]{\begin{propos}{#1}\sl }
\newcommand{\eprop}{\end{propos}}

\newenvironment{proposs}[1]%
    { \begin{PropS} \label{P:#1}}%
    { \end{PropS} }
\newcommand{\props}[1]{\begin{proposs}{#1}\sl }
\newcommand{\eprops}{\end{proposs}}

\newenvironment{corol}[1]%
    { \begin{Cor} \label{C:#1}}%
    { \end{Cor} }
\newcommand{\cor}[1]{\begin{corol}{#1} \sl }
\newcommand{\ecor}{\end{corol}}

\newenvironment{corols}[1]%
    { \begin{CorS} \label{C:#1}}
    { \end{CorS} }
\newcommand{\cors}[1]{\begin{corols}{#1} \sl }
\newcommand{\ecors}{\end{corols}}

\newenvironment{defeni}[1]%
    { \begin{Def} \label{D:#1}}%
    { \end{Def} }
\newcommand{\defe}[1]{\begin{defeni}{#1} }
\newcommand{\edefe}{\end{defeni}}

\newenvironment{defenis}[1]%
    { \begin{DefS} \label{D:#1}}%
    { \end{DefS} }
\newcommand{\defes}[1]{\begin{defenis}{#1} }
\newcommand{\edefes}{\end{defenis}}

\newenvironment{remark}[1]%
    { \begin{Rem} \label{R:#1}}%
    { \end{Rem} }
\newcommand{\rem}[1]{\begin{remark}{#1}}
\newcommand{\erem}{\end{remark}}

\newenvironment{remarks}[1]%
    { \begin{RemS} \label{R:#1}}%
    { \end{RemS} }
\newcommand{\rems}[1]{\begin{remarks}{#1}}
\newcommand{\erems}{\end{remarks}}

\newenvironment{conjec}[1]%
    { \begin{Conj} \label{Co:#1}}
    { \end{Conj} }
\renewcommand{\conj}[1]{\begin{conjec}{#1} \sl }
\newcommand{\econj}{\end{conjec}}

\newenvironment{example}[1]%
    { \begin{Ex} \label{Exx:#1}}%
    { \end{Ex} }
\newcommand{\ex}[1]{\begin{example}{#1} }
\newcommand{\eex}{\end{example}}

\newenvironment{examples}[1]%
    { \begin{ExS} \label{Exx:#1}}
    { \end{ExS} }
\newcommand{\exs}[1]{\begin{examples}{#1} }
\newcommand{\eexs}{\end{examples}}

\newcommand{\prf}{ \begin{proof} }
\newcommand{\epr}{ \end{proof} }

\newcommand{\eu}{\begin{center}
\begin{minipage}{5cm}
\includegraphics[width=4.8cm]{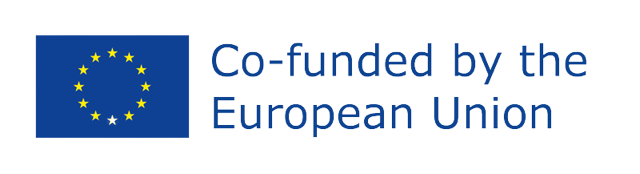}
\end{minipage}
\begin{minipage}{9.5cm}
\footnotesize Co-Funded by the European Union. Views and opinions expressed are however those of the author(s) only and do not necessarily reflect those of the European Union. Neither the European Union nor the granting authority can be held responsible for them.
\end{minipage}
\end{center}}

\newcommand\Z{\mathbb{Z}}
\newcommand\Q{\mathbb{Q}}

\newcommand\calF{\mathcal{F}}
\newcommand\calO{\mathcal{O}}
\newcommand\calI{\mathcal{I}}

\mathchardef\mhyphen="2D

\newcommand\pmeninf{p^{-\infty}}
\DeclareMathOperator{\ten}{\otimes}
\DeclareMathOperator{\im}{Im}

\newcommand{\Gm}{\mathbb{G}_m}
\newcommand{\Ga}{\mathbb{G}_a}
\newcommand{\Fp}{\mathbb{F}_p}
\newcommand{\Zp}{\mathbb{Z}_p}
\newcommand{\Qp}{\mathbb{Q}_p}
\newcommand{\mup}{\mu_{p}}                                              
\newcommand{\mupn}{\mu_{p^n}}                                           

\newcommand{\onto}[2][]{\xrightarrow[#1]{#2}\mathrel{\mkern-14mu}\rightarrow}           
\newcommand\xto[2]{\xrightarrow[#1]{#2}}

\DeclareMathOperator{\id}{id}                
\DeclareMathOperator{\Hom}{Hom}                                      
\DeclareMathOperator{\coker}{coker}

\DeclareMathOperator{\Art}{Art}
\DeclareMathOperator{\D}{D}                                              
\DeclareMathOperator{\Tot}{Tot}                                          
\DeclareMathOperator{\Cone}{Cone}                                        
\DeclareMathOperator{\Ab}{\textit{Ab}}                                   
\DeclareMathOperator{\Rgam}{\text{R}\Gamma}                              
\DeclareMathOperator{\Spec}{Spec}                                        
\newcommand{\W}{W}                                                       
\newcommand{\K}{K}
\DeclareMathOperator{\Rlim}{Rlim}                                        
\DeclareMathOperator{\Rder}{\text{R}}                                    
\DeclareMathOperator{\Rflat}{R\Gamma_{\textnormal{fppf}}}                
\DeclareMathOperator{\RdR}{R\Gamma_{\textnormal{dR}}}                    
\DeclareMathOperator{\Pic}{Pic}                                          
\DeclareMathOperator{\Br}{Br}                                            

\newcommand\Hflat[1]{H^{#1}_{\textnormal{fppf}}}                         
\newcommand\Het[1]{H^{#1}_{\textnormal{\'{e}t}}}                         

\newcommand\HdR[1]{H^{#1}_{\textnormal{dR}}}                             

\DeclareMathOperator{\Acr}{\mathbb{A}_{cris}}                            
\newcommand\perf[1]{#1^{\flat}}                                          
\newcommand\coperf[1]{#1_{\textnormal{perf}}}                            
\DeclareMathOperator{\Rcris}{R\Gamma_{cris}}                             
\newcommand{\Nyg}{\mathrm{F}^1_\mathrm{N}}                                      
\newcommand\Nygm[1]{\mathrm{F}_{\mathrm{N}}^{1,#1}}                                      
\newcommand\Hcris[1]{H^{#1}_{\textnormal{cris}}}                         
\newcommand\fg[2]{\Phi^{#1}(#2,\Gm)}
\newcommand\ffg[2]{\Phi^{#1}_{\mathrm{fl}}(#2,\Gm)}
\newcommand\sitofppf[1]{\left(#1\right)_{\mathrm{fppf}}}

\usepackage[style=french]{csquotes}
\newcommand\slogan[1]{
\begin{center}
\begin{minipage}{0.9\textwidth}
\begin{center}
    \enquote{\emph{#1}}
\end{center}
\end{minipage}
\end{center}
\vspace{.3em}}

\begin{document}

\title{Formal smoothness of the Artin--Mazur formal groups}
\author{Livia Grammatica\footnote{Email: \textit{livia.grammatica@math.unistra.fr}\\
Institut de Recherche Math\'{e}matique Avanc\'{e}e (IRMA), Universit\'{e} de Strasbourg, 7 rue Ren\'{e} Descartes, 67000 Strasbourg, France.}}  
\date{\today}

\maketitle
\begin{abstract}
Let $X$ be a smooth proper variety over an algebraically closed field of positive characteristic $p$. We find cohomological conditions for the Artin--Mazur formal group functors $\fg{i}{X}$ to be formally smooth. We show that if all crystalline cohomology groups of $X$ are torsion-free (e.g. if $X$ is an abelian variety) then all of the $\fg{i}{X}$ are representable and formally smooth. We then identify a necessary condition for formal smoothness, which we use to give examples, for any $d\ge2$, of varieties $X$ for which $\fg{i}{X}$ is formally smooth when $i<d$, whereas $\fg{d}{X}$ is not. The constructions are inspired by Igusa's surface with non-smooth Picard scheme. Finally, we give a condition equivalent to formal smoothness in terms of Serre's Witt vector cohomology. The strategy relies on the notion of $C$-smoothness - where $C$ is the group algebra of $\Qp/\Zp$ - which is a condition that detects when a formal group is formally smooth, and on the use of the Nygaard filtration to relate fppf cohomology to crystalline cohomology.
\end{abstract}
\tableofcontents

\sec{}{Introduction}

Let $k$ be an algebraically closed field of positive characteristic $p$. For a smooth and proper $k$-variety $X$, Artin and Mazur \cite{artinmazur} define a family of formal group functors $\fg{i}{X}$ by the rule
\begin{equation*}
    \fg{i}{X}(R)=\ker\left(\Het{i}(X_R,\Gm)\to\Het{i}(X_{R_{\mathrm{red}}},\Gm)\right),
\end{equation*}
for any Artinian $k$-algebra $R$ (where $X_R$ stands for the base change $X\times_k R)$. The functor $\fg{1}{X}$ is the formal completion at the origin of the Picard scheme of $X$, while $\fg{2}{X}$ is usually called the formal Brauer group of $X$, denoted by $\widehat{\mathrm{Br}}(X)$. Note that $\fg{1}{X}$ is always representable, but this may not hold for $\fg{i}{X}$ when $i\ge2$. For example, by \cite[Proposition 10.11]{braggolsson}, $\widehat{\mathrm{Br}}(X)$ is representable if and only if $\Pic_{X/k}$ is smooth. We are interested in the following question of Artin--Mazur \cite[p.104]{artinmazur}: \slogan{assuming that $\fg{i}{X}$ is representable, when is it formally smooth?}
In the special case $i=1$, this amounts to asking when $\Pic_{X/k}$ is smooth. One of the goals of this paper is to construct, in characteristic $2$ and for any $d\ge2$, a variety $Z$ for which $\fg{d}{Z}$ is representable but non-smooth. To be precise, we prove the following. 

\th{esempiintro}\hspace{-.2em}\textnormal{(\Cref{P:esempietto}, \Cref{T:esempio})}
Suppose $p=2$ and let $d\ge2$ be an integer. There exists a smooth proper $k$-variety $Z$ such that 
\begin{enumerate}
    \item $\fg{i}{Z}$ is representable for $i\le d$,
    \item $\fg{i}{Z}$ is formally smooth for $i<d$,
    \item $\fg{d}{Z}$ is not formally smooth.
\end{enumerate}
When $d=2$ there is an explicit $Z$ satisfying these conditions.
\eth

We construct $Z$ as follows: let $E$ be an ordinary elliptic curve over $k$, with its non-zero $2$-torsion point $a$. If $Y$ is a smooth proper variety over $k$, equipped with an involution $\tau:Y\to Y$, we get an involution $\sigma$ of $E\times Y$ mapping $(x,y)$ to $(x+a,\tau(y))$. Then the quotient $Z=E\times Y/\langle\sigma\rangle$ is a smooth and proper variety -- the reader may notice a similarity with Igusa's surface \cite{igusa}, see \Cref{SubS:esempi}. 

To prove that, under suitable hypothesis on $(Y,\tau)$, $Z$ has the desired properties, we establish two simple cohomological criteria for the formal smoothness of the Artin--Mazur formal groups. Let us go back to $X$ being an arbitrary smooth proper $k$-variety. The first result is an easily verifiable sufficient condition for $\fg{i}{X}$ to be formally smooth, in terms of the torsion of crystalline cohomology. 

\th{sufficienteintro}\hspace{-.2em}\textnormal{(\Cref{T:sufficiente})} If $\fg{i}{X}$ is representable, and $\Hcris{i+1}(X/\W(k))$ is torsionfree, then $\fg{i}{X}$ is formally smooth.
\eth

Combining this with a representability result of Artin--Mazur (\Cref{P:rappartinmazur}), it follows that for an abelian variety $A$ all functors $\fg{i}{A}$ are representable and formally smooth (\Cref{C:varab}).

The next result gives a necessary condition in terms of the torsion of the $p$-adic étale cohomology groups $\Het{i}(X,\Zp)$.

\th{necessariaintro}\hspace{-.2em}\textnormal{(\Cref{T:necessaria})} Suppose that $\fg{i}{X}$ is representable and that $\Het{i+1}(X,\Zp)_{\text{tors}}\ne0$. Then $\fg{i}{X}$ is not formally smooth.
\eth

Using the identification of $\Het{i}(X,\Zp)$ with $\Hcris{i}(X,\Zp)^{F=1}$ (\cite[II.5.2]{illusiedrw}) this can also be seen as a criterion on the torsion of crystalline cohomology. 

Via the Hochschild-Serre spectral sequence for crystalline and étale cohomology, we can use these results to find conditions on the cohomology of $Y$, and on the involution $\tau$, for $Z$ to satisfy the conditions of \Cref{T:esempiintro}. The final step is then to actually construct pairs $(Y,\tau)$ with the required cohomology. This is done in \Cref{SubS:esempi}.

We also present one result which is independent from the examples discussed above. It identifies a condition equivalent to the formal smoothness of $\fg{i}{X}$, in terms of Serre's Witt vector cohomology \cite{serremexico}. Recall that
\begin{equation*}
    H^{d}(X,\mathcal{W})=\varprojlim_n H^{d}(X,\mathcal{W}_n),
\end{equation*}
for $d\ge0$, where $\mathcal{W}_n$ is the truncated Witt vector group scheme, and the transition maps are induced by the restriction maps $\mathcal{W}_{n+1}\to \mathcal{W}_n$. Thus the groups $H^d(X,\mathcal{W})$
are endowed with a Frobenius endomorphism $F$ and a Verschiebung $V$, and they sit in a long exact sequence
\begin{equation*}
\dots\to H^d(X,\mathcal{W})\xto{}{V} H^d(X,\mathcal{W})\to H^{d}(X,\calO_X)\to H^{d+1}(X,\mathcal{W})\to\dots,
\end{equation*}
see \Cref{SubS:pruffafinale} for details.

\th{equivalenteintro}\hspace{-.2em}\textnormal{(\Cref{T:equivalente})} Suppose that $\fg{i}{X}$ is representable. The following are equivalent.

$(1)$ the formal group $\fg{i}{X}$ is formally smooth.

$(2)$ the map $H^{i}(X,\mathcal{W})\to H^{i}(X,\calO_X)$ is surjective.

$(3)$ the Verschiebung $V$ acting on $H^{i+1}(X,\mathcal{W})$ is injective.

\eth

Although our proofs of the three criteria are original, the statements can also be deduced from existing results which rely on completely different techniques. \Cref{T:equivalenteintro} is found in work of \cite{ekedahl}, see also \cite[Section 6.1-6.2]{yuan} and \cite[Sections 10-12]{braggolsson}. Ekedahl \cite[Proposition III.8.1]{ekedahl} proves that there is a short exact sequence of $\W(k)$-modules
\begin{equation*}
0\to\mathbb{D}(\ffg{i}{X}_{\mathrm{inf}})\to H^{i+1}(X,\mathcal{W})\to\mathbb{D}(\ffg{i+1}{X}_{\mathrm{sm}})\to0
\end{equation*}
where $\mathbb{D}$ is the covariant Dieudonné module functor. Here $\ffg{i}{X}$ is the fppf sheafification of $\fg{i}{X}$, and agrees with $\fg{i}{X}$ when the latter is representable, see \Cref{SubS:gruppiformali} for details. Since $V$ is injective on the right-hand group and nilpotent on the left-hand group, one has \Cref{T:equivalenteintro}. Ekedahl's paper makes use of deep properties of the de Rham-Witt complex.

\Cref{T:sufficienteintro} follows from \Cref{T:equivalenteintro} by an advanced result in the theory of the de Rham-Witt complex (the ``survie du c{\oe}ur'', \cite[II.3.4]{illusieraynaud}), and the condition of \Cref{T:necessariaintro} is implied by condition $(3)$, as can be seen by an elementary Artin-Schreier argument. However, we wanted to highlight these two results because Witt vector cohomology is not as well-behaved, nor as easy to compute, as crystalline or étale $p$-adic cohomology (see e.g. \cite[II.7]{illusiedrw}). Moreover, we do not use the de Rham-Witt complex in any of our arguments, which instead rely on the systematic use of quasisyntomic descent and the Nygaard filtration.

\vspace{.5cm}
\noindent\textbf{Strategy of proof and outline.} Let us sketch how we obtain our cohomological criteria. In \Cref{SubS:liscezza} we introduce the notion of $C$-smoothness. Let $C=k[\Qp/\Zp]$, and let $f_C:C\to C$ be the $k$-algebra map induced by multiplication by $p$ on $\Qp/\Zp$. If $\calF$ is an abelian sheaf on $\sitofppf{k}$, or a formal group over $k$, it makes sense to consider the endomorphism $\calF(f_C)$ of $\calF(\Spec(C))$. We say that $\calF$ is $C$\textit{-smooth} if $\calF(f_C)$ is surjective (\Cref{D:defcliscezza}). 

This definition is a useful tool to study the formal smoothness of formal groups: we prove in \Cref{P:cliscezza2} that if $G$ is a formal group of finite type, then $G$ is formally smooth if and only if $G$ is $C$-smooth. In \Cref{SubS:gruppiformali} we apply this to the Artin--Mazur groups, and via the Kummer sequence one shows that the formal smoothness of the $\fg{i}{X}$ (when it is representable) is equivalent to the $C$-smoothness of the fppf sheaf $\varinjlim_n\Rder^if_{*}\mu_{p^n,X}$, see \Cref{P:equivalenzacliscezza}. We further show in \Cref{C:criteriofinale} that this is equivalent to the $f_C$ acting surjectively on $\Hflat{i}(X\times\Spec(C),\Qp/\Zp(1))$. The latter sits in a short exact sequence
\begin{equation*}
0\to\Hflat{i}(X_C,\Zp(1))\otimes\Qp/\Zp\to\Hflat{i}(X_C,\Qp/\Zp(1))\to\Hflat{i+1}(X_C,\Zp(1))_{\text{tors}}\to0,
\end{equation*}
which splits non-canonically.

The relation with crystalline cohomology comes from the fact that fppf cohomology with $\Zp(1)$ coefficients can be described in terms of crystalline cohomology and of the Nygaard filtration - we explain this in \Cref{SubS:cristallina}. This allows us to translate our previous smoothness criterion into one involving crystalline invariants of $X$. In \Cref{SubS:criterio} we prove \Cref{T:sufficienteintro} and \Cref{T:necessariaintro}. In \Cref{SubS:esempi} we use these results to construct the examples announced earlier in the introduction, and the proof of \Cref{T:equivalenteintro} is given in the final section of the paper.

\vspace{.5cm}
\noindent\textbf{Quasisyntomic descent.} This technique is used without much introduction in \Cref{S:pruffa}, to prove many results about crystalline cohomology, therefore we briefly describe it here. A remarkable property of crystalline cohomology is that if $f:Y\to Z$ is a quasisyntomic cover of quasisyntomic $k$-schemes (\cite[Definition 4.10]{bms2}), then the map
\begin{equation*}   
\xymatrix{\Rcris(Z/\W(k))\ar[r]&\Tot\Bigl(\Rcris(Y/\W(k))\ar[r]<1.5pt>\ar[r]<-1.5pt> & \Rcris(Y\times_Z Y/\W(k))\ar[r]\ar[r]<-3pt>\ar[r]<3pt> & \dots\Bigr)}
\end{equation*}
in $\D_{\ge0}(\Zp)$ is an isomorphism. The same property holds for other crystalline invariants such as the Nygaard filtration, and also fppf cohomology with value in a finite type affine group scheme - see \cite{bhattlurie} for an exhaustive account, or \cite{articolouno} for a more elementary reference. Descent makes it possible to reduce questions on the crystalline cohomology of $Z$, to the same question for each $n$-fold fibre product of $Y$ over $Z$. A simple usecase is reducing a problem for separated smooth schemes to a problem for affine smooth schemes. 

Descent is especially useful if the crystalline cohomology of $Y$ of $Y\times_ZY,$ etc. is easy to describe. If $Z=\Spec(A)$ is affine, the projection $Y=\coperf{Z}=\varprojlim_{F_Z}Z\to Z$ is a quasisyntomic cover and every $n$-fold fibre product of $Y$ over $Z$ is the spectrum of a qrsp $k$-algebra. The point is that the crystalline cohomology of $\Spec(B)$ where $B$ is qrsp is isomorphic to $\Acr(B)[0]$, and is well-understood (see \cite[8.2]{bms2}).

The most relevant case to us is when $Z=\Spec(A)$ is an affine smooth $k$-scheme. Then $\coperf{A},\coperf{A}\otimes_A\coperf{A}$ are what we call elementary qrsp $k$-algebras. The functor $\Acr(-)$ (and the Nygaard filtration) for these algebras can be described very explicitly. In this paper, we will use this to reduce many statements about the crystalline cohomology of smooth schemes, to statements about the (semi-)linear algebra on the rings $\Acr(B)$ for $B$ elementary qrsp. More details and references are given in \Cref{SubSubS:eqrsp}.

\vspace{.5cm}
\noindent\textbf{Relation to existing results.} A number of results regarding the smoothness of $\Pic_{X/k}$ can be found in the literature. Igusa \cite{igusa} first constructed a variety $X$ with non-smooth Picard scheme, over a field of characteristic $2$. In \cite{serremexico} the author introduces Witt vector cohomology and studies its main properties. He constructs along the way varieties with non-smooth Picard group in any characteristic, for which Hodge symmetry fails. Building on this, Mumford \cite{mumford} shows that for a smooth proper $X$, $\Pic_{X/k}$ is smooth if and only if the map
\begin{equation*}
    H^1(X,\mathcal{W})\to H^1(X,\calO_X)
\end{equation*}
is surjective, which is the case $i=1$ of \Cref{T:equivalenteintro}. Berthelot and Nygaard \cite{nygaard} prove that if $\Hcris{2}(X/\W(k))$ is torsion-free then $\Pic_{X/k}$ is smooth, which is the $i=1$ case of \Cref{T:sufficienteintro}.

The Artin--Mazur formal groups are higher-dimensional versions of the Picard scheme. They were first studied in \cite{artinssk3} as an interesting invariant of $K3$ surfaces in positive characteristic. When $X$ is a $K3$ surface, the formal Brauer group $\widehat{\Br}(X)$ is a smooth $1$-dimensional formal group classified by its height $h$. The height of $X$ determines its crystalline cohomology \cite[II.7.2]{illusiedrw} and, when $h=\infty$, the variety $X$ is supersingular \cite{artinssk3}, a notion with no analogue in characteristic zero. More recently \cite{cyss} these considerations have been extended to study the geometry of Calabi-Yau varieties $X$ in positive characteristic: if $X$ is $n$-dimensional, $\fg{n}{X}$ is a smooth $1$-dimensional formal group. 

Artin--Mazur \cite{artinmazur} computed the tangent space of $\fg{i}{X}$ as $H^i(X,\calO_X)$, so if this group is $0$ then $\fg{i}{X}$ is formally smooth for trivial reasons. They also show that if $H^{i+1}(X,\calO_X)=0$ then $\fg{i}{X}$ is formally smooth. Both of these results are a consequence of \Cref{T:equivalenteintro}, using the equivalence between condition $(2)$, and Serre's Bockstein homomorphisms \cite{serremexico} being all equal to $0$.

The most general results concerning Artin--Mazur formal groups are proved by Ekedahl \cite{ekedahl}. As mentioned above, he proves a refinement of \Cref{T:equivalenteintro}, and much more. Yuan Yang kindly explained to the author that many of the results presented here can be obtained using deep results in the theory of the de Rham-Witt complex found in Illusie--Raynaud's and Ekedahl's work.

\vspace{.5cm}
\noindent\textbf{Acknowledgements.} I wish to express heartfelt thanks to my advisor Marco D'Addezio, who introduced me to this problem and suggested the strategy to tackle it. This work would not have been possible without his encouragements and many enriching discussions with him. I benefited from helpful discussions with Emiliano Ambrosi, with Yuan Yang, who taught me a different approach to the Artin--Mazur groups, and special thanks are due to Yuya Matsumoto, who helped me with one of the examples of \Cref{SubS:esempi}. Thanks are also due to Alexander Petrov and Emanuel Reinecke for useful comments. Most of this work was done while attending the Bernoulli program "Arithmetic geometry of K3 surfaces" at EPFL, I am very grateful for the invitation. This project has received funding from the European Union’s Horizon Europe research
and innovation programme under the Marie Skłodowska-Curie grant agreement n° 101126554.

\vspace{.5cm}
\noindent\textbf{Notation.} We fix the prime number $p$ and an algebraically closed field $k$ of characteristic $p$ throughout the paper, with the exception of \Cref{SubS:esempi} where we take $p=2$. If $A$ is a $k$-algebra denote by $\W(A)$ (resp. $\W_n(A)$) the (truncated) Witt vectors with values in $A$. Write $\Art_k$ for the category of Artinian $k$-algebras, and if $Y$ is a $k$-scheme write $\sitofppf{Y}$ for the big fppf site of $Y$. We will denote by $F$ the absolute Frobenius of $Y$, and also every map induced by $F$ by functoriality. This should not cause confusion as we only ever work with one scheme at a time. 

Whenever we write "derived category" we mean the derived $\infty$-category in the sense of Lurie. This robust framework allows us to define maps and objects by quasisyntomic descent without having to worry about categorical issues.

\sec{preliminari}{Preliminaries on formal groups}

After a brief reminder on formal groups, in \Cref{SubS:liscezza} we give the definition $C$-smoothness and $C$-étaleness, and we show that for finite type formal groups they agree with the usual  conditions of smoothness and étaleness (\Cref{P:cliscezza2}). In \Cref{SubS:gruppiformali} we recall some properties of the Artin--Mazur formal groups and use the results of the previous section to show that formal smoothness of Artin--Mazur formal groups is equivalent to a condition on fppf cohomology groups with $\Zp(1)$-coefficients.

\ssec{liscezza}{Formal groups and $C$-smoothness}

Let $C=k[\Qp/\Zp]$, which is isomorphic to $k[y^{\pmeninf}]/(y-1)$. We take $x=y-1$ and thus identify $C$ with $k[x^{\pmeninf}]/(x)$. The notion of $C$-smoothness is designed to determine if a formal group (or group scheme) $G$ is formally smooth by looking at its $C$-points. Recall that a formal group functor is a functor $G:\Art_k\to\Ab$, and a formal group is a formal group functor $G$ which is representable, meaning that there is a profinite $k$-algebra $A=\varprojlim A_n$ and an isomorphism of functors
\[
G\simeq\Hom_{\mathrm{cont}}(A,-).
\]
We say that $G$ is connected if $A$ is local, étale if $A$ is a product of copies of $k$, of finite type if $A$ is topologically finitely generated as a $k$-algebra. For a formal group $G$ there is a canonical short exact sequence
\[
0\to G^0\to G\to G_{\text{ét}}\to0
\]
where $G^0$ is connected and $G_{\text{\'{e}t}}$ is étale.

A formal group functor $G$ can be extended to the category of ind-objects of $\Art_k$ by setting
\begin{equation}
    G(R)=\varinjlim G(R_n)\text{ for }R=(R_n,f_n)
\end{equation}
and similarly for maps. Then writing $G(C)$ makes sense: we can write $C=\varinjlim A_m$ where $A_m=k[x]/(x^{p^m})$ and the transition maps send $x$ to $x^p$.

The $k$-linear morphism 
\begin{equation}
f_C:C\to C \quad\quad x^{1/p^i}\mapsto x^{1/p^{i-1}}
\end{equation}
will be ubiquitous in this note. As a morphism of ind-Artin algebras, it is the colimit of the maps $p_m:A_{m+1}\to A_m$ mapping $x$ to $x$.

\defe{defcliscezza}
Let $\calF$ be a presheaf on $\sitofppf{k}$ or a formal group functor over $k$. We say that
\begin{enumerate}
    \item $\calF$ is $C${-smooth} if $\calF(f_C)$ is surjective,
    \item $\calF$ is $C${-etale} if $\calF(f_C)$ is bijective.
\end{enumerate}

\edefe

Denote by $\mathfrak{G}_p$ be the category of abelian finite type affine group schemes over $k$, which are killed by some power of $p$. Any group $G\in\mathfrak{G}_p$ is a sheaf over $\sitofppf{k}$ via the Yoneda embedding. It turns out that for $G\in\mathfrak{G}_p$, being $C$-smooth or $C$-étale is the same as being smooth or étale in the usual sense. We need the following technical lemma.

\lem{vanishingcoom}

If $G\in\mathfrak{G}_p$ then $\Hflat{i}(\Spec(C),G)=0$ for $i>0$. 

\elem

\prf
Any $G\in\mathfrak{G}_p$ is a successive extension of the groups $\Ga,\alpha_p,\Z/p$ and $\mup$, and if $L$ is any of these groups we have $\Hflat{i}(\Spec(C),L)=0$ for $i>0$. This is clear for $\Ga$. Using the short exact sequences
\begin{equation}
    0\to\alpha_p\to\Ga\xto{}{F}\Ga\to0
\end{equation}
\begin{equation}
    0\to\Z/p\to\Ga\xto{}{F-1}\Ga\to0
\end{equation}
we find 
\begin{equation}\label{alphap}
    \Hflat{1}(\Spec(C),\alpha_p)=\coker\left(F_C:C\to C\right),
\end{equation}
\begin{equation}\label{zetap}
    \Hflat{1}(\Spec(C),\Z/p)=\coker\left(F_C-1:C\to C\right),
\end{equation}
where $F_C:C\to C$ is the absolute Frobenius. The group \eqref{alphap} is zero because $C$ is semiperfect. As for the group \eqref{zetap}: if $a\in k\subseteq C$, $a$ is in the image of $1-F_C$ because $k$ is algebraically closed. If $a\in\ker\left(C\onto[]{}k\right)$, the sum $x=a+f_C(a)+f_C^2(a)+\dots$ is finite, and $(1-f_C)(x)=a$. Therefore this group is also $0$. For $\mup$ this is less elementary, it is proved in \cite[Proposition 7.2.5]{bhattlurie}. 

The lemma now follows by dévissage.
\epr

\rem{}
    It can be shown by dévissage that if $G\in\mathfrak{G}_p$ is not the trivial group then $G(C)\ne0$. Therefore $G\mapsto G(C)$ is a faithful embedding of $\mathfrak{G}_p$ into the category of abelian groups.
\erem

\prop{cliscezza1}

Let $G\in\mathfrak{G}_p$. The following hold:
\begin{enumerate}
    \item $G$ is $C$-smooth if and only if $G$ is smooth.
    \item $G$ is $C$-étale if and only if $G$ is étale, i.e. a discrete $p$-group.
\end{enumerate}

\eprop

\prf

If $0\to H\to G\to G/H\to0$ is a short exact sequence in $\mathfrak{G}_p$, by \Cref{L:vanishingcoom} we have a commutative diagram 
\begin{equation}\label{fCgruppi}
\begin{tikzcd}
	0 & {H(C)} & {G(C)} & {G/H(C)} & 0 \\
	0 & {H(C)} & {G(C)} & {G/H(C)} & 0
	\arrow[from=1-1, to=1-2]
	\arrow[from=1-2, to=1-3]
	\arrow["{H(f_C)}", from=1-2, to=2-2]
	\arrow[from=1-3, to=1-4]
	\arrow["{G(f_C)}", from=1-3, to=2-3]
	\arrow[from=1-4, to=1-5]
	\arrow["{G/H(f_C)}", from=1-4, to=2-4]
	\arrow[from=2-1, to=2-2]
	\arrow[from=2-2, to=2-3]
	\arrow[from=2-3, to=2-4]
	\arrow[from=2-4, to=2-5]
\end{tikzcd}
\end{equation}
with exact rows. So if the outermost vertical map are surjective, the middle vertical map is also surjective. And if the middle vertical map is surjective so is the rightmost vertical map.

If $G$ is smooth, then $G$ is a successive extension of $\Ga$ and $\Z/p$. It is easy to check that $f_C$ is surjective on $\Ga(C)$ and $\Z/p(C)$, so using diagram (\ref{fCgruppi}) inductively we find that $G(f_C):G(C)\to G(C)$ is surjective.

Suppose now that $G$ is infinitesimal and non-trivial. We prove by induction on the length of $G$ that $G(f_C):G(C)\to G(C)$ is not surjective. If $G$ is of length $1$ then $G=\alpha_p$ or $G=\mup$, and $G(f_C)=0$. If $G$ is of length $>1$ it has a non-trivial subgroup $H$. The group $G/H$ is infinitesimal, non-trivial and of smaller length, so $G/H(f_C)$ is not surjective. Thus $G(f_C)$ can not be surjective.

Finally, for general $G$, consider the short exact sequence
\[
0\to G_{\mathrm{red}}\to G\to G_{\mathrm{inf}}\to0
\]
where $G_{\mathrm{red}}$ is smooth and $G_{\mathrm{inf}}$ is infinitesimal. If $f_C:G(C)\to G(C)$ is surjective so is $f_C:G_{\mathrm{inf}}(C)\to G_{\mathrm{inf}}(C)$. By the previous paragraph we must have $G_{\mathrm{inf}}=0$ i.e. $G$ is smooth. This proves (1), and the non-trivial part of (2) is proved similarly. \epr

Extending \Cref{P:cliscezza1} to formal groups requires an additional small argument.

\prop{cliscezza2}
Let $G$ be a formal group whose identity component $G^0$ is of finite type. The following hold: 
\begin{enumerate}
    \item $G$ is $C$-smooth if and only if $G$ is formally smooth.
    \item $G$ is $C$-étale if and only if $G$ is étale.
\end{enumerate}
\eprop

\prf

To prove (1) we may suppose that $G$ is connected. Suppose that $G$ is formally smooth. Then each of the maps $G(p_m):G(A_{m+1})\to G(A_m)$ is surjective. Passing to the colimit we see that $G(f_C)$ is surjective.

Viceversa, suppose that $G(f_C)$ is surjective. Consider the short exact sequence 
\[
0\to G_{\mathrm{red}}\to G\to G_{\mathrm{inf}}\to0
\]
where $G_{\mathrm{red}}$ is formally smooth and $G_{\mathrm{inf}}$ is a finite connected formal group, i.e. a finite infinitesimal group scheme. Then $G\to G_{\mathrm{inf}}$ is formally smooth, so for all $m$ the map $G(A_m)\to G_{\mathrm{inf}}(A_m)$ is surjective. Passing to the limit we find that $G(C)\to G_{\mathrm{inf}}(C)$ is surjective. So we have a commutative diagram
\[\begin{tikzcd}
	{G(C)} & {G_{\mathrm{inf}}(C)} \\
	{G(C)} & {G_{\mathrm{inf}}(C)}
	\arrow[from=1-1, to=1-2]
	\arrow["{G(f_C)}", from=1-1, to=2-1]
	\arrow["{G_{\mathrm{inf}}(f_C)}", from=1-2, to=2-2]
	\arrow[from=2-1, to=2-2]
\end{tikzcd}\]
where the horizontal maps and the leftmost vertical map are surjective. It follows that $G_{\mathrm{inf}}(f_C)$ is surjective, and \Cref{P:cliscezza1} now implies that $G_{\mathrm{inf}}=0$, i.e. $G$ is formally smooth.

The non-trivial part of point (2) follows from \Cref{P:cliscezza1}. 
\epr

\ssec{gruppiformali}{The Artin--Mazur formal groups} 

For a smooth projective variety $X$, Artin--Mazur \cite{artinmazur} define a family of formal group functors $\fg{i}{X}$ by the rule
\begin{equation*}
    \fg{i}{X}(R)=\ker\left(\Het{i}(X_R,\Gm)\to\Het{i}(X_{R_{\mathrm{red}}},\Gm)\right),
\end{equation*}
where $X_R$ is shorthand for $X\times\Spec(R)$. For $i=1$ this is the formal group associated to $\Pic_{X/k}$ and it is always representable. As previously noted in the introduction, for $i\ge2$ the functor $\fg{i}{X}$ may fail to be representable. On the other hand, Raynaud \cite{raynaud} and Bragg-Olsson \cite{braggolsson} prove that if we let $\ffg{i}{X}$ denote the fppf sheafification of $\fg{i}{X}$, the group functor $\ffg{i}{X}$ is always representable. The tangent space of $\ffg{i}{X}$ is finite-dimensional (see \Cref{P:compconn} below) so $\ffg{i}{X}$ is a connected formal group of finite type. Thus we prefer to work with $\ffg{i}{X}$ throughout, while keeping in mind that when $\fg{i}{X}$ is representable it coincides with $\ffg{i}{X}$. Bragg--Olsson give the following criterion for $\fg{i}{X}$ to be representable.

\prop{rappartinmazur} (\cite[Proposition 10.11]{braggolsson})
The following are equivalent.

(1) $\fg{i}{X}$ is representable.

(2) $\ffg{i-1}{X}$ is formally smooth.

\eprop

Via the Kummer sequence one can relate $\ffg{i}{X}$ to a sheaf which looks perhaps more natural. Let $f:X\to\Spec(k)$ be the structure map. Avatars of the following theorem go back to Milne \cite{milnemupn} but the statement presented here is due to Bragg-Olsson.

\th{rappmupn}(\cite[Theorem 1.3]{braggolsson}) The fppf sheaf $\Rder^if_*\mupn{}$ is representable by a group in $\mathfrak{G}_p$.
\eth

Given a group $G\in\mathfrak{G}_p$ its formal completion at the origin, denoted $\widehat{G}$, is the connected formal group defined by the rule
\begin{equation*}
    R\mapsto\ker\left(G(R)\to G(R_{\mathrm{red}})\right)
\end{equation*}
for any Artin algebra $R$. The tangent space of $\widehat{G}$ is naturally isomorphic to the tangent space at the origin of $G$. If $G^i_{n}$ denotes the formal completion of $\Rder^if_*\mupn{}$ at the origin, we thus obtain a formal group $G^i=\varinjlim_n G^i_{n}$.

\prop{compconn}

The formal groups $G^i$ and $\ffg{i}{X}$ are naturally isomorphic, and the tangent space of $\ffg{i}{X}$ is a quotient of $H^i(X,\calO_X)$. 

\eprop

\prf The first part is proved via the Kummer sequence, see \cite[Proposition 10.7]{braggolsson}. For the second part, use \cite[Proposition 10.10]{braggolsson}, and the fact that the tangent space to $\fg{i}{X}$ is $H^i(X,\calO_X)$, \cite[Corollary 2.4]{artinmazur}.
\epr

Let's apply \Cref{P:cliscezza2} in this setting.

\prop{equivalenzacliscezza}
The following are equivalent:
\begin{enumerate}
    \item $\ffg{i}{X}$ is formally smooth
    \item the sheaf $\varinjlim_n\Rder^if_{*}\mupn{}$ is $C$-smooth.
\end{enumerate}
\eprop

\rem{}
    Let $(\calF_i)_{i\in\calI}$ be a filtered diagram of sheaves on $\sitofppf{k}$. For a $k$-scheme $T$ the natural map
    \begin{equation*}
        \varphi_T:\varinjlim_i(\calF_i(T))\to(\varinjlim_i\calF_i)(T)
    \end{equation*}
    is usually not bijective. But if $T$ is quasi-compact and quasi-separated (e.g. affine), it satisfies condition (4) of \cite[0738]{stacksproject}, so $\varphi_T$ is an isomorphism. We will use this tacitly in following calculations.
\erem

\prf

By \Cref{P:compconn}, $\ffg{i}{X}$ is formally smooth if and only if $G^i$ is formally smooth, which by \Cref{P:cliscezza2} is equivalent to $G^i$ being $C$-smooth. Now,
\begin{align*}
G^i(\Spec(C))&=\varinjlim_m\varinjlim_n\ker\left(\Rder^if_{*}\mupn{}(\Spec(A_m))\to\Rder^if_{*}\mupn{}(\Spec(k))\right)\\
&=\varinjlim_n\varinjlim_m\ker\left(\Rder^if_{*}\mupn{}(\Spec(A_m))\to\Rder^if_{*}\mupn{}(\Spec(k))\right)\\
&=\ker\left(\varinjlim_n\varinjlim_m\Rder^if_{*}\mupn{}(\Spec(A_m))\to\varinjlim_n\varinjlim_m\Rder^if_{*}\mupn{}(\Spec(k))\right)\\
&=\ker\left(\varinjlim_n\Rder^if_{*}\mupn{}(\Spec(C))\to\varinjlim_n\Rder^if_{*}\mupn{}(\Spec(k))\right),
\end{align*}
where in the last step we used the equality
\begin{equation*}\label{conseguenzafinitario}
    \Rder^if_{*}\mupn{}(\Spec(C))=\varinjlim_m\Rder^if_{*}\mupn{}(\Spec(A_m)).
\end{equation*}
This holds because, in the terminology of \cite{marcobrauer}, $\Rder^if_{*}\mupn{}$ is finitary, see Lemma 3.3 for a proof. 

Therefore, since $k\to C$ has a retraction, we have proved that there is a direct sum decomposition
\begin{equation*}
    \varinjlim_n\Rder^if_{*}\mupn{}(\Spec(C))=G^i(\Spec(C))\oplus\varinjlim_n\Rder^if_{*}\mupn{}(\Spec(k))
\end{equation*}
Since $f_C$ acts identically on the right-hand summand, the $C$-smoothness of $G^i$ is equivalent to (2), and finally (1) and (2) are equivalent. \epr

The point of this maneuver is that $\varinjlim_n\Rder^if_{*}\mupn{}(C)$ can be identified with an fppf cohomology group. Then, via the comparison of fppf cohomology with crystalline cohomology, we will establish a crystalline criterion for $f_C$ to be surjective. To prove this identification we need some auxiliary results on the cohomology of fppf sheaves on $\Spec(C)$.

\prop{}

Let $\pi:X_C\to\Spec(C)$ be the projection on the first factor. 
\begin{enumerate}
    \item[(1)] The sheaves $\Rder^if_*\mupn{}|_C$ and $\Rder^i\pi_*\mupn{}$ on the big fppf site of $\Spec(C)$ are isomorphic.
    \item[(2)] The map $\Hflat{i}(X_C,\mupn{})\to\Hflat{0}(\Spec(C),\Rder^i\pi_*\mupn{})\simeq\Rder^if_*\mupn{}(C)$ coming from the Leray spectral sequence of $\pi$ is an isomorphism.
\end{enumerate}

\eprop

\prf

Point (1) is formal: by definition, both $\Rder^if_*\mupn{}|_C$ and $\Rder^i\pi_*\mupn{}$ can be described as the sheafification of
\begin{equation*}
V\mapsto\Hflat{i}(X\times_kV,\mupn),
\end{equation*}
so they are isomorphic.

For (2) it suffices to show that $\Hflat{j}(\Spec(C),\Rder^i\pi_*\mupn{})$ vanishes for all $j>0$. Using that $\Rder^if_*\mupn{}$ is representable by some group in $\mathfrak{G}_p$ (\Cref{T:rappmupn}), this is a consequence of (1) and \Cref{L:vanishingcoom}. \epr

We get the following criterion for formal smoothness in terms of fppf cohomology.

\cor{ccriteriofinale}

The formal group $\ffg{i}{X}$ is formally smooth if and only if the action of $\id_X\times f_C$ on 
\begin{equation*}
\Hflat{i}(X_C,\Qp/\Zp(1))=\varinjlim_n\Hflat{i}(X_C,\mupn)
\end{equation*}
is surjective. \qed

\ecor

The group $\Hflat{i}(X_C,\Qp/\Zp(1))$ sits in a short exact sequence
\begin{equation}\label{succqsuz}
0\to\Hflat{i}(X_C,\Zp(1))\otimes\Qp/\Zp\to\Hflat{i}(X_C,\Qp/\Zp(1))\to\Hflat{i+1}(X_C,\Zp(1))_{\text{tors}}\to0
\end{equation}
of $\Zp$-modules. In the next section, we describe the outermost groups in terms of crystalline cohomology, which allows us to study the action of $\id_X\times f_C$ in more detail.

\sec{pruffa}{Relation with crystalline cohomology and main results}

The previous section culminated in \Cref{C:ccriteriofinale}, where we related the formal smoothness of the Artin--Mazur groups to a condition on the groups $\Hflat{i}(X_C,\Zp(1))$. The goal of this section is to translate this condition into one involving the crystalline cohomology of $X$. To this end, we use the description of fppf cohomology of $X_C$ as the fibre of a map between the first piece of the Nygaard filtration of $X_C$ and the crystalline cohomology of $X_C$ (\Cref{T:comparisonbl}). Thus, in \Cref{SubS:cristallina} we study in detail the crystalline cohomology and the Nygaard filtration of $X_C$, using the modern framework of \cite{bms2} and \cite{bhattlurie} which relies importantly on quasisyntomic descent. This allows us to prove in \Cref{SubS:criterio} and \Cref{SubS:pruffafinale} the results stated in the introduction. In \Cref{SubS:esempi}, we construct, for any $d$, varieties $X$ for which $\fg{d}{X}$ is representable but not formally smooth, while $\fg{i}{X}$ is representable and formally smooth for $i<d$. Throughout \Cref{SubS:cristallina}, \Cref{SubS:criterio} and \Cref{SubS:pruffafinale} we fix a smooth proper $k$-variety $X$. 

\ssec{cristallina}{Fppf cohomology and crystalline cohomology}

For a $k$-scheme $Y$, denote by $\Rcris(Y/\W(k))$ its crystalline cohomology as an object of $\D(\W(k))$. We write
\begin{equation*}
    F:\Rcris(Y/\W(k))\to\Rcris(Y/\W(k))
\end{equation*}
for the action of the absolute Frobenius, considered as a $\sigma$-linear map in $\D(\W(k))$ or as a linear map in $\D(\Zp)$. The fibre of the augmentation map 
\begin{equation}\label{defnygaard}
\Rcris(Y/\W(k))\to\Rgam(Y,\calO_Y)
\end{equation}
is called the first piece of the Nygaard filtration on $\Rcris(Y/\W(k))$, and denoted by $\Nyg{}\Rcris(Y/\Zp)$. Thus we have a long exact sequence of $\W(k)$-modules
\begin{equation}\label{sel2}
	\cdots \to {H^{i-1}(Y,\calO_Y)} \to {\Nyg\Hcris{i}(Y/\W(k))} \to {\Hcris{i}(Y/\W(k))} \to {H^i(Y,\calO_Y)} \to \cdots,
\end{equation}
and since $\Rgam(Y,\calO_Y)$ is killed by $p$, after inverting $p$ the complexes $\Nyg{}\Rcris(Y/\W(k))$ and $\Rcris(Y/\W(k))$ become isomorphic. We denote by $\iota$ the natural map from the Nygaard filtration to crystalline cohomology.

Now let $Y$ be quasisyntomic - see \cite[Definition 4.10]{bms2} for a definition, for our purposes we only need that smooth $k$-schemes, $\Spec(C)$ and products of these are quasisyntomic. An important feature of the Nygaard filtration is the existence of a divided Frobenius $F/p:\Nyg{}\Rcris(Y/\W(k))\to\Rcris(Y/\W(k))$ making the diagram
\[\begin{tikzcd}
	{\Nyg{}\Rcris(Y/\W(k))} & {\Rcris(Y/\W(k))} \\
	{\Rcris(Y/\W(k))} & {\Rcris(Y/\W(k))}
	\arrow["{F/p}", from=1-1, to=1-2]
	\arrow["\iota", from=1-1, to=2-1]
	\arrow["p", from=1-2, to=2-2]
	\arrow["F", from=2-1, to=2-2]
\end{tikzcd}\]
commute (\cite[5.3.3]{bhattlurie}). We can now state the main result comparing fppf cohomology with crystalline cohomology via the Nygaard filtration. It goes back to \cite{fontainemessing} but see also \cite[7.3.5]{bhattlurie} for a proof with quasisyntomic descent.

\th{comparisonbl}
There is an exact triangle
\begin{equation}
    \Rflat(Y,\Zp(1))\to\Nyg{}\Rcris(Y/\W(k))\xto{}{F/p-\iota}\Rcris(Y/\W(k))
\end{equation}
in $\D(\Zp)$.
\eth

Thus to understand the fppf cohomology of $X_C$ we need a good grasp of crystalline cohomology, the Nygaard filtration, and $F/p-\iota$ for the schemes $X,C,$ and $X_C$. The rest of this section is devoted to understanding this, but before we move on we make one last remark on the map $F/p-\iota$. As of now, to alleviate notation we will write $\Hcris{i}(Y)$ instead of $\Hcris{i}(Y/\W(k))$, and similarly for the Nygaard filtration.

From \eqref{sel2} we can extract the short exact sequences
\begin{equation}\label{nyginquadrata}
    0\to S^i_Y\to H^i(Y,\calO_Y)\to V^i_Y\to0
\end{equation}
where $S^i_Y$ comes from $\Hcris{i}(Y)$ and $V^i_Y$ is a submodule of $\Nyg{}\Hcris{i+1}(Y)$. Studying the action of $F/p$ on $V^i_Y$ is an important step toward understanding the action on $\Nyg{}\Hcris{i+1}(Y)$ whole.

Consider the map $\Rgam(Y,\calO_Y)\to\Nyg{}\Rcris(Y)[1]$ coming from \eqref{defnygaard}. Since $p$ acts as zero on $\Rgam(Y,\calO_Y)$, it factors through the fibre of multiplication by $p$ on $\Nyg{}\Rcris(Y)[1]$ (here it is important that we work with $\infty$-categories to get a well-defined factorization). In other words, we have a commutative diagram
\begin{equation}\label{diagrammafbar}
\begin{tikzcd}
	{\Rgam(Y,\calO_Y)} \\
	{\Nyg{}\Rcris(Y)\otimes^L\Z/p} & {\Nyg{}\Rcris(Y)[1]} & {\Nyg{}\Rcris(Y)[1]} \\
	{\Rcris(Y)\otimes^L\Z/p} & {\Rcris(Y)[1]} & {\Rcris(Y)[1]}
	\arrow[dashed, from=1-1, to=2-1]
	\arrow[from=1-1, to=2-2]
	\arrow[from=2-1, to=2-2]
	\arrow["{F/p}"', from=2-1, to=3-1]
	\arrow["{p}", from=2-2, to=2-3]
	\arrow["{F/p}"', from=2-2, to=3-2]
	\arrow["{F/p}"', from=2-3, to=3-3]
	\arrow[from=3-1, to=3-2]
	\arrow["{p}", from=3-2, to=3-3]
\end{tikzcd}\end{equation}
Let $\overline{F}:\Rgam(Y,\calO_Y)\to\Rcris(Y)\otimes^L\Z/p$ be the composition of the left-hand vertical maps. Then understanding $F/p-1$ on $V^i_Y$ amounts to understanding the composition
\begin{equation}\label{fbar}
H^{i-1}(Y,\calO_Y)\xto{}{\overline{F}} H^{i-1}(\Rcris(Y)\otimes^L\Z/p)\to\Hcris{i}(Y),
\end{equation}
so we will need to study the behaviour of $\overline{F}$.

\lem{}
The diagram
\begin{equation}\label{triangoloffbar}
\begin{tikzcd}
	{\Rcris(Y)\otimes^L\Z/p} && {\Rcris(Y)\otimes^L\Z/p} \\
	& {\Rgam(Y,\calO_Y)}
	\arrow["F", from=1-1, to=1-3]
	\arrow[from=1-1, to=2-2]
	\arrow["{\overline{F}}"', from=2-2, to=1-3]
\end{tikzcd}\end{equation}
commutes, where the unnamed map is the augmentation map \eqref{defnygaard}.
\elem
\prf
We only sketch the proof, similar arguments will be detailed later on. By quasisyntomic descent we may reduce to $Y=\Spec(R)$ for some qrsp algebra $R$. Then one must show that 
\[
\Acr(R)/p\to R\xto{}{\overline{F}}\Acr(R)/p
\]
is equal to the absolute Frobenius. This holds because $\Acr(R)/p$ is the divided power envelope of $\ker\left(\perf{R}\onto{}R\right)$, \cite[F.7]{bhattlurie}.
\epr

\rem{remarka}
We see easily that $\overline{F}:R\to\Acr(R)/p$ is a ring homomorphism. Therefore $\overline{F}$ is compatible with the Künneth isomorphism.
\erem

\sssec{}{Cohomology of $X$}

If $X$ is proper, $\Rcris(X/\W(k))$ is a perfect object of $\D(\W(k))$ (\cite[Tag 07MX]{stacksproject}), in particular $\Hcris{i}(X/\W(k))$ is a $\W(k)$-module of finite type. The rest of the results presented in this short discussion only require $X$ to be smooth, so the reader may assume this. The crystalline-de Rham comparison theorem \cite[Tag 07MI]{stacksproject} gives an isomorphism 
\begin{equation*}
    \Rcris(X/\W(k))\otimes^L\Z/p\simeq\RdR(X/k)
\end{equation*}
and thus a short exact sequence
\begin{equation}\label{crisderham}
    0\to\Hcris{i}(X/\W(k))/p\to\HdR{i}(X/k)\to\Hcris{i+1}(X/\W(k))[p]\to0
\end{equation}
for all $i$. Before describing $\overline{F}$ we recall a construction in de Rham cohomology: for a smooth $k$-scheme $Y$ there is a map of complexes $\calO_Y[0]\xto{}{}\Omega^{\bullet}_{Y}$ given by the absolute Frobenius in degree $0$. It induces a Frobenius-linear morphism
\begin{equation*}
    F_{\mathrm{dR}}:\Rgam(Y,\calO_Y)\to\RdR(Y/k).
\end{equation*}

\lem{fsupcoerente}

The map $\overline{F}:H^i(X,\calO_X)\to\Hcris{i+1}(X)$ of \eqref{fbar} is equal to the composition
\begin{equation}\label{comp2}
H^i(X,\calO_X)\xto{}{F_{\mathrm{dR}}}\HdR{i}(X/k)\to\Hcris{i+1}(X)[p]\subseteq\Hcris{i+1}(X)
\end{equation}
\elem

\prf
We must show that the two maps $\overline{F},F_{\mathrm{dR}}$ are equal. By Zariski descent it is enough to prove it for $Y=\Spec(R)$ affine and smooth. Then the only nonzero cohomology groups are in degree $0$. By design the map $H^0(F_{\mathrm{dR}}):R\to\HdR{0}(Y/k)=R^p$ is the $p$-the power map. On the other hand, taking $H^0$ in \eqref{triangoloffbar} we get a commutative diagram
\begin{equation*}
    \begin{tikzcd}
	{R^p} && {R^p} \\
	& {R}
	\arrow["^\wedge p", from=1-1, to=1-3]
	\arrow[from=1-1, to=2-2]
	\arrow["{H^0(\overline{F})}"', from=2-2, to=1-3]
    \end{tikzcd}
\end{equation*}
and since $R$ is reduced, the ring map $H^0(\overline{F})$ is uniquely determined and equal to $H^0(F_{\mathrm{dR}})$. \epr

\sssec{eqrsp}{Cohomology of $\Spec(C)$.} The crystalline cohomology of quasiregular semimperfect schemes such as $\Spec(C)$ is well-documented. We recall the basic facts and refer to \cite{drinfeldacris} and \cite{articolouno} for more details.

If $A$ is a semiperfect $k$-algebra, $\Rcris(\Spec(A))$ is isomorphic to $\Acr(A)[0]$. From \eqref{defnygaard} we see that the Nygaard filtration is isomorphic to
\[
\Nyg\Acr(A)=\ker\left(\Acr(A)\onto{}A\right).
\]
concentrated in degree $0$. If furthermore  
\begin{equation}\label{eqrsp}
    A\simeq B[x_1^{\pmeninf},\dots,x_n^{\pmeninf}]/(x_1,\dots,x_n)
\end{equation}
for some perfect $k$-algebra $B$ we have an explicit description of $\Acr(A)$ and of the Nygaard filtration. Keeping the terminology of \cite{articolouno}, we say $A$ is an elementary quasiregular semiperfect $k$-algebra, eqrsp for short. We introduce some notation.

If $\alpha\in\Z_+[1/p]$, denote by $(\alpha!)_p$ the largest power of $p$ dividing ($\lfloor\alpha\rfloor!)$. If $\alpha=(\alpha_1,\dots,\alpha_n)$ is a vector of elements of $\Z_+[1/p]$ we will write $(\alpha!)_p$ for the product of the $(\alpha_i!)_p$. If furthermore $x_1,\dots,x_n$ is a set of indeterminates we will be concerned with the monomials 
\begin{equation*}
    x^{<\alpha>}=\frac{x_1^{\alpha_1}\dots x_n^{\alpha_n}}{(\alpha_1!)_p\dots(\alpha_n!)_p}
\end{equation*}
which we abbreviate as $x^{<\alpha>}$. Thus, for example, we have identities such as
\begin{equation*}
    x^{<\alpha>}x^{<\beta>}=x^{<\alpha+\beta>}\frac{((\alpha+\beta)!)_p}{(\alpha!)_p(\beta!)_p}.
\end{equation*}

\lem{descrizioneacris}
(1) We can identify $\Acr(A)$ with the ring of power series of the form 
\begin{equation*}
    \sum_{\alpha\in\Z_+[1/p]^n}b_{\alpha}x^{<\alpha>},\quad b_{\alpha}\in \W(B)
\end{equation*} 
such that for every $n>0$ the set $\{\alpha\text{ s.t. }p^n\text{ does not divide }b_{\alpha}\}$ is finite. Thus $\Acr(A)$ is flat over $\W(B)$.

(2) The Frobenius $F$ acts as 
\[
F\left(\sum_{\alpha\in\Z_+[1/p]^n}b_{\alpha}x^{<\alpha>}\right)=\sum_{\alpha\in\Z_+[1/p]}F(b_{\alpha})p^{\lfloor\alpha\rfloor} x^{<p\alpha>}.
\] 
(3) The ideal $\Nyg\Acr(C)$ is the set of power series $\sum_{\alpha\in\Z_+[1/p]^n}b_{\alpha}x^{<\alpha>}$ such that $p$ divides $b_{\alpha}$ whenever $\alpha_i<1$ for all $i$.\qed
\elem

We will use this notation whenever we encounter eqrsp rings in the rest of the article. Specialising to $C$, we get an explicit description of $\Acr(C)$ and of $\Nyg{}\Acr(C)$. We will need to consider the action of $f_C$ on $\Acr(C)$, which we denote by $a\mapsto f_C^{\mathrm{cris}}(a)$. Since the absolute Frobenius of $C$ is the composition of $f_C$ with the Frobenius of $k$, it follows from \Cref{L:descrizioneacris} that
\begin{equation*}
    f_C^{\mathrm{cris}}\left(b_{\alpha}x^{<\alpha>}\right)=b_{\alpha}p^{\lfloor\alpha\rfloor}x^{<p\alpha>}.
\end{equation*}

As for $\overline{F}$, following diagram \eqref{diagrammafbar} we find that it is Frobenius-linear and maps $x^{\alpha}\in C$ to $x^{<p\alpha>}\in\Acr(C)$. 

We finish with a technical result which will be used later on.

\lem{conucleoacris}

Let $s,r\ge0$ be coprime integers with $s>0$, and set 
\[
\mathbb{A}_{\mathrm{cris}}^{r,s}(C)=\{a\in\Acr(C)|p^{r-s}F^s(a)\in\Acr(C)\}
\]
The map $\mathbb{A}_{\mathrm{cris}}^{r,s}(C)\to\Acr(C)$ defined by $a\mapsto p^{r-s}F^s(a)-a$ is surjective.

\elem

\prf

It is sufficient to check surjectivity mod $p$, so let $M:\mathbb{A}_{\mathrm{cris}}^{r,s}(C)\ten\Fp\to\Acr(C)\ten\Fp$ be the reduction mod $p$ of the map in the statement. From \Cref{L:descrizioneacris} we see that $\Acr(C)\ten\Fp$ is the free $k$-module with basis all $x^{<\alpha>}$. Therefore it is enough to show that if $b\in k$ and $\alpha\in\Z_+[1/p]$ then $a=bx^{<\alpha>}$ is in the image of $M$.

If $\sum_{i=0}^{s-1}\lfloor\alpha p^{i}\rfloor>s-r$ then $M(-a)=a$ and we are done. If $\sum_{i=0}^{s-1}\lfloor\alpha p^{i}\rfloor=s-r$, then 
\begin{equation*}
    a+M(a)=b^{p^s}x^{<p^s\alpha>},
\end{equation*}
which we have just shown to be in the image of $M$. So $a$ is also in the image of $M$. Suppose that $\sum_{i=0}^{s-1}\lfloor\alpha p^{i}\rfloor<s-r$ and call $\ell$ the (positive) difference of these two integers. Then 
\begin{equation*}
    c=p^{\ell}b^{1/p^s}x^{<\alpha/p^s>}
\end{equation*} 
is nonzero in $\mathbb{A}_{\mathrm{cris}}^{r,s}(C)\otimes\Fp$ and $M(c)=a$.\epr

\sssec{}{Cohomology of $X_C$.} The crystalline cohomology of $X_C$ is computed via the crystalline Künneth formula.

\prop{kunnethacriscomp}
There is a canonical quasi-isomorphism
\[
\Rcris(X_C)\simeq\Rcris(X)\widehat{\otimes}^{L}_{\W(k)}\Acr(C)
\]
where Frobenius on the left corresponds to $F\otimes F$ on the right.
\eprop

\prf

This is a very special case of \cite[4.1.8]{bhattlurie}. We sketch a direct proof based on quasi-syntomic descent.

Note that there is a canonical map from the right-hand complex to the left-hand complex. Passing to an appropriate quasisyntomic cover of $X$, we reduce to showing that if $A$ is eqrsp, as in \eqref{eqrsp}, the canonical map 
\[
\Acr(A)\widehat{\otimes}_{\W(k)}\Acr(C)\to\Acr(A\otimes_kC)
\]
is an isomorphism. As both sides are $p$-complete it is enough to show that
\begin{equation*}
    \Acr(A)/p^n\otimes_{\W_n(k)}\Acr(C)/p^n\to\Acr(A\otimes_kC)/p^n
\end{equation*}
is an isomorphism. This is straightforward using the explicit descriptions of $\Acr(A)$ and $\Acr(C)$ given in \Cref{L:descrizioneacris}.
\epr

\cor{kunnethacriscoom}

For all $i\ge0$ there is a canonical isomorphism
\[
\Hcris{i}(X_C)\simeq\Hcris{i}(X)\otimes_{\W(k)}\Acr(C),
\]
where Frobenius on the left corresponds to $F\otimes F$ on the right.

\ecor

\prf

Use \Cref{P:kunnethacriscomp}, the flatness of $\Acr(C)$ as a $\W(k)$-module, and that $\Rcris(X)$ is a perfect complex of $\W(k)$-modules.\epr

The Künneth formula for coherent cohomology likewise gives isomorphisms
\begin{equation*}
    H^i(X_C,\calO_{X_C})\simeq H^{i}(X,\calO_X)\otimes_kC
\end{equation*}
for all $i$, hence the long exact sequence
\begin{equation*}
    \dots\to H^{i-1}(X,\calO_X)\otimes_kC\to\Nyg{}\Hcris{i}(X_C)\to\Hcris{i}(X)\otimes_{\W(k)}\Acr(C)\to\cdots.
\end{equation*}
In the notation of \eqref{nyginquadrata} we see that $V^i_{X_C}=V^i_X\otimes_kC$. The following lemma explains the action of $F/p$ on this group.

\lem{fbarprodotto}
The map
\[
\overline{F}:H^i(X,\calO_X)\otimes_kC\to\left(\Hcris{i+1}(X)\otimes_{\W(k)}\Acr(C)\right)[p]=\Hcris{i+1}(X)[p]\otimes_k\Acr(C)/p
\]
is equal to $\overline{F}\otimes\overline{F}$.
\elem

\prf
This follows from \Cref{R:remarka}.
\epr

For now, let $Y,Z$ be any two quasisyntomic schemes. The composition
\begin{align*}
    \Nyg{}\Rcris(Y)\widehat{\otimes}^L_{\W(k)}\Rcris(Z)\to\Rcris(Y)&\widehat{\otimes}^L_{\W(k)}\Rcris(Z)\\
    &\simeq\Rcris(Y\times Z)\to\Rgam(Y\times Z,\calO_{Y\times Z})
\end{align*}
is zero, thereby producing a map 
\begin{equation*}
    i_{Y,Z}:\Nyg{}\Rcris(Y)\widehat{\otimes}^L_{\W(k)}\Rcris(Z)\to\Nyg{}\Rcris(Y\times Z).
\end{equation*}
We have the following commutativity property.
\lem{}
The diagram
\begin{equation}\label{nygaardtensore}
\begin{tikzcd}
	{\Nyg{}\Rcris(Y)\widehat{\otimes}^L_{\W(k)}\Rcris(Z)} & {\Rcris(Y)\widehat{\otimes}^L_{\W(k)}\Rcris(Z)} \\
	{\Nyg{}\Rcris(Y\times Z)} & {\Rcris(Y\times Z)}
	\arrow["{F/p\otimes F}", from=1-1, to=1-2]
	\arrow["i_{Y,Z}", from=1-1, to=2-1]
	\arrow["\simeq", from=1-2, to=2-2]
	\arrow["{F/p}", from=2-1, to=2-2]
\end{tikzcd}
\end{equation}
commutes.
\elem

\prf
By quasisyntomic descent it is enough to check this on affine quasiregular semiperfect algebras, which is a straightforward verification.
\epr

If we specialize to $X$ and $C$ we have, for any $i\ge0$, the two maps
\begin{gather*}
    i_{X,C}:\Nyg{}\Hcris{i}(X)\otimes_{\W(k)}\Acr(C)\to\Nyg{}\Hcris{i}(X_C),\\
    i_{C,X}:\Hcris{i}(X)\otimes_{\W(k)}\Nyg{}\Acr(C)\to\Nyg{}\Hcris{i}(X_C),
\end{gather*}
which we can use to describe completely $\Nyg{}\Hcris{i}(X_C)$.

\lem{lemnygaardprodotto} 
The following hold.

$(1)$ The cokernel of $i_{X,C}:\Nyg{}\Hcris{i}(X)\otimes_{\W(k)}\Acr(C)\to\Nyg{}\Hcris{i}(X_C)$ is killed by $\id_X\times f_C$. 

$(2)$ Take $a_1,\dots,a_n\in\Hcris{i}(X)$ which form a basis when projected on $S^i_X$, in the notation of \eqref{nyginquadrata}. Using the notation of \Cref{L:descrizioneacris}, any element of $\Nyg{}\Hcris{i}(X_C)$ can be written uniquely in the form
\begin{equation}\label{formaelt}
i_{X,C}\left(\sum_{\alpha<1}a_{\alpha}\otimes x^{\alpha}+\sum_{\alpha\ge1}a_{\alpha}\otimes x^{<\alpha>}\right)+i_{C,X}\left(\sum_{\alpha\ge1}b_{\alpha}\otimes x^{<\alpha>}\right),
\end{equation}
where the sums converge $p$-adically, $a_{\alpha}\in\Nyg{}\Hcris{i}(X)$ for all $\alpha$, and is determined mod $V^{i-1}_X$ when $\alpha\ge1$, and each $b_{\alpha}$ is of the form $[\lambda_1]a_1+\dots+[\lambda_n]a_n$ for some $\lambda_i\in k$.

$(3)$ The map $F/p$ maps an element of the form \eqref{formaelt} to
\begin{equation}\label{formulafsup}
    \sum_{\alpha}\left(p^{\lfloor\alpha\rfloor}(F/p)(a_\alpha)\otimes x^{<p\alpha>}\right)+\sum_{\alpha\ge1}\left(p^{\lfloor\alpha\rfloor-1}F(b_\alpha)\otimes x^{<p\alpha>}\right)
\end{equation}
in $\Hcris{i}(X)\otimes\Acr(C)$. The map $\id_X\times f_C$ maps an element of the form \eqref{formaelt} to
\begin{equation}\label{formulafC}
i_{X,C}\left(\sum_{\alpha<1}a_{\alpha}\otimes  x^{<p\alpha>}+\sum_{\alpha\ge1}a_{\alpha}\otimes  x^{<p\alpha>}\right)+i_{C,X}\left(\sum_{\alpha\ge1}b_{\alpha}\otimes x^{<p\alpha>}\right).
\end{equation}
in $\Nyg{}\Hcris{i}(X_C)$.
\elem

\prf

$(1)$ Consider the diagram
\begin{equation}\label{diagrammanecessaria}
\begin{tikzcd}
	{H^{i-1}(X,\calO_X)\otimes C} & {\Nyg{}\Hcris{i}(X_C)} & {\Hcris{i}(X_C)} \\
	{H^{i-1}(X,\calO_X)\otimes\Acr(C)} & {\Nyg{}\Hcris{i}(X)\otimes\Acr(C)} & {\Hcris{i}(X)\otimes\Acr(C)}
	\arrow[from=1-1, to=1-2]
	\arrow[from=1-2, to=1-3]
	\arrow["g", two heads, from=2-1, to=1-1]
	\arrow[from=2-1, to=2-2]
	\arrow["j", from=2-2, to=1-2]
	\arrow[from=2-2, to=2-3]
	\arrow["\simeq"', from=2-3, to=1-3]
\end{tikzcd}
\end{equation}
where the top row is \eqref{sel2} for $X_C$ and the bottom row comes from tensoring \eqref{sel2} for $X$ with $\Acr(C)$ (tensors products are taken over $\W(k)$). Thus both rows are exaxt. 

We know that $\Nyg{}\Hcris{i}(X_C)$ surjects onto 
\begin{equation*}
    T:=\ker(\Hcris{i}(X)\otimes\Acr(C)\to H^{i+1}(X,\calO_X)\otimes C).
\end{equation*}
A simple diagram chase (using the fact that $g$ is surjective) shows that $\coker(j)$ is isomorphic to
\begin{align*}
T/\im\left(\Nyg{}\Hcris{i}(X)\otimes\Acr(C)\to\Hcris{i}(X)\otimes\Acr(C)\right)\\
\simeq\im\left(S^{i+1}_X\otimes\Nyg{}\Acr(C)\to S^{i+1}_X\otimes\Acr(C)\right),
\end{align*}
hence $\id_X\times f_C$ acts as zero on this group, because $S^{i+1}_X$ is $p$-torsion.

$(2)$ Chasing diagram \eqref{diagrammanecessaria} we see that the kernel of $j$ is equal to 
\begin{equation*}
    \im\left(V^{i-1}_X\otimes\Nyg{}\Acr(C)\to V^{i-1}_X\otimes\Acr(C)\right)\subseteq\Nyg{}\Hcris{i}(X)\otimes\Acr(C).
\end{equation*}
Moreover, it is easy to check that any element of $\coker j$, which we described in $(1)$, lifts uniquely to an element 
\begin{equation*}
    \sum_{\alpha\ge1}b_{\alpha}x^{<\alpha>}
\end{equation*}
such as in the statement of the lemma. The conclusion follows easily from the explicit description of $\Acr(C)$ and $\Nyg{}\Acr(C)$ of \Cref{L:descrizioneacris}.

$(3)$ Follows from \eqref{nygaardtensore}, and the fact that $i_{X,C}$ and $i_{C,X}$ are equivariant with respect to $\id_X\times f_C$.
\epr

\ssec{criterio}{Proof of \Cref{T:sufficienteintro} and \Cref{T:necessariaintro}}

Recall that we want to understand the action of $f_C$ on $\Hflat{i}(X_C,\Zp(1))$. In the previous section we related these groups to other groups, namely $\Nyg{}\Hcris{i}(X_C)$ and $\Hcris{i}(X_C)$, which we described in detail. From this study we derive in \Cref{C:criteriofinale} a necessary and sufficient condition for $\ffg{i}{X}$ to be formally smooth in terms of crystalline cohomology and the Nygaard filtration of $X_C$. This in turn allows us to prove \Cref{T:sufficienteintro} and \Cref{T:necessariaintro} stated in the introduction.

\lem{sintomicaeqrsprazionale}
The cokernel of 
\begin{equation}\label{cokernelfsup}
F/p-\iota:\Nyg{}\Hcris{i}(X_C)\to\Hcris{i}(X_C)
\end{equation}
is torsion for all $i$. Thus $\Hflat{i}(X_C,\Qp(1))$ is isomorphic to
\begin{equation*}
    \Hcris{i}(X_C)[1/p]^{F=p}=\left(\Hcris{i}(X)[1/p]\otimes_{\K}\Acr(C)[1/p]\right)^{F\otimes F=p}
\end{equation*}
for all $i$.
\elem

\prf

We need to show that 
\[
\frac{F\otimes F}{p}-1:\Hcris{i}(X)[1/p]\otimes_{\K}\Acr(C)[1/p]\to\Hcris{i}(X)[1/p]\otimes_{\K}\Acr(C)[1/p]
\]
is surjective for all $i$.

If $\Hcris{i}(X)[1/p]=M_1\oplus\dots\oplus M_n$ is a decomposition in simple isocrystals, $F\otimes F/p-1$ preserves the subspaces $M_i\otimes_{\W(k)}\Acr(C)[1/p]$, so it is enough to check that
\[
F\otimes F/p-1:M\otimes_{\K}\Acr(C)[1/p]\to M\otimes_{\K}\Acr(C)[1/p]
\]
is surjective when $M$ is a simple isocrystal of slope $\lambda=r/s$. Choose a basis $x_1,\dots,x_s$ of $M$ such that 
\[
F(x_1)=x_2,\dots,F(x_{s-1})=x_s, F(x_s)=p^rx_1.
\]
Then we have to check that the map
\begin{align*}
&F\otimes F/p-1:\Acr(C)[1/p]^s\to\Acr(C)[1/p]^s\\
&(a_1,\dots,a_s)\mapsto\left(p^{r-1}F(a_s)-a_1,F(a_1)/p-a_2,\dots,F(a_{s-1})/p-a_s\right)
\end{align*}
is surjective, i.e. if $(b_1,\dots,b_s)\in\Acr(C)[1/p]$, we look for $a=(a_1,\dots,a_s)$ such that $(F\otimes F/p-1)(a)=b$. Solving this system, we find that the only obstruction is solving the equation
\[
p^{r-s}F^s(a_s)-a_s=b_s+F(b_{s-1})p^{-1}+\dots+F^{s-1}(b_1)p^{-s+1}
\]
By \Cref{L:conucleoacris} a solution exists, so we are done. \epr

One can prove along the same lines that the cokernel of \eqref{cokernelfsup} has finite $p$-exponent.

\prop{suriettivosinistra}
The map $\id_X\times f_C$ acts surjectively on $\Hflat{i}(X_C,\Qp(1))$. Therefore it acts surjectively on $\Hflat{i}(X_C,\Zp(1))\otimes\Qp/\Zp$.
\eprop

\prf

There is a map $s:X_C\to X_C$ such that
\[
s\circ(\id_X\times f_C)=(\id_X\times f_C)\circ s=F_{X_C}
\]
where $F_{X_C}$ is the absolute Frobenius of $X_C$. \Cref{L:sintomicaeqrsprazionale} shows that $F_{X_C}$ acts as multiplication by $p$ on $\Hflat{i}(X_C,\Qp(1))$, which is bijective. Therefore $\mathrm{id}_X\times f_C$ also induces a bijection on $\Hflat{i}(X_C,\Qp(1))$.\epr

Now \eqref{succqsuz} tells us that $\ffg{i}{X}$ is formally smooth if and only if $\id_X\times f_C$ is surjective on $\Hflat{i+1}(X_C,\Zp(1))_{\mathrm{tors}}$.
This group sits in a short exact sequence
\[
0\to\coker\left(F/p-\iota\right)\to\Hflat{i+1}(X_C,\Zp(1))_{\mathrm{tors}}\to\ker\left(F/p-\iota\right)_{\text{tors}}\to0,
\]
as a consequence of \Cref{L:sintomicaeqrsprazionale}.

\prop{suriettivodestra}

The map $\id\times f_C$ acts surjectively on $\coker(F/p-\iota)$.

\eprop

\prf

The cokernel of
\[
F/p-\iota:\Nyg{}\Hcris{i}(X_C)\to\Hcris{i}(X_C)=\Hcris{i}(X)\otimes_{\W(k)}\Acr(C)
\]
is torsion (\Cref{L:sintomicaeqrsprazionale}) and generated as a $\W(k)$-module by classes of the form $a=m\otimes x^{<\alpha>}$ where $m\in\Hcris{i}(X)$. Recall that the action of $\id_X\times f_C$ is given by
\begin{equation*}
    a\mapsto m\otimes f_C^{\mathrm{cris}}\left(x^{<\alpha>}\right)=m\otimes p^{\lfloor\alpha\rfloor} x^{<p\alpha>}
\end{equation*}
If $\alpha< p$ we have $a=(\id_X\times f_C)(m\otimes b)$ for $b=x^{<\alpha/p>}$.

If $\alpha\ge 1$, then $a\in\Hcris{i}(X)\otimes\Nyg{}\Acr(C)$, thus $a$ is in the same class as $a_1=a+(F/p-1)(i_{C,X}(a))$ in $\coker\left(F/p-1\right)$. Now $a_1=F(m)\otimes p^{\lfloor\alpha\rfloor-1} x^{<p\alpha>}$, which also lies in $\Hcris{i}(X)\otimes_{\W(k)}\Nyg{}\Acr(C)$, so we can repeat the argument and see that $a$ is in the same class as
\[
a_{\ell}:=F^{\ell}(m)\otimes\left(p^{\lfloor\alpha\rfloor+\dots+\lfloor p^{\ell-1}\alpha\rfloor-\ell}x^{<p^l\alpha>}\right)
\]
for every ${\ell}\ge1$. Choose ${\ell}\gg 0$ such that $\lfloor\alpha\rfloor+\dots+\lfloor p^{{\ell}-1}\alpha\rfloor-{\ell}>p^{{\ell}-1}\alpha$. As before, we can select $b\in\Acr(C)$ such that $(\id_X\times f_C)(F^{\ell}(m)\otimes b)=a_{\ell}$, and we are done. \epr

\cor{criteriofinale}
The fomal group $\ffg{i}{X}$ is formally smooth if and only if $\id_X\times f_C$ acting on $\ker\left(F/p-\iota\right)_{\mathrm{tors}}$ is surjective.\qed
\ecor

\th{sufficiente}

If $\Hcris{i+1}(X)$ is torsionfree, $\ffg{i}{X}$ is formally smooth.

\eth

\prf
Under this hypothesis, we have
\begin{equation*}
    \Nyg{}\Hcris{i+1}(X_C)_{\mathrm{tors}}=V^{i+1}_X\otimes C
\end{equation*}
in the notation of \Cref{SubS:cristallina}. Moreover the map $F/p-1$ restricted to this subgroup is zero, because the target is torsion-free. Therefore,
\begin{equation*}
    \ker\left(F/p-\iota\right)_{\mathrm{tors}}=V^{i+1}_X\otimes C,
\end{equation*}
and $\id_X\times f_C$ acts surjectively on this group because $f_C$ acts surjectively on $C$. By \Cref{C:criteriofinale}, $\ffg{i}{X}$ is formally smooth.
\epr

\cor{varab}

If $A$ is an abelian variety, the Artin--Mazur functors $\fg{i}{A}$ are all representable and formally smooth.

\ecor
\prf
Combine \Cref{P:rappartinmazur} and \Cref{T:sufficiente}.
\epr

To conclude the section, we give is our sufficient condition for the formal group $\ffg{i}{X}$ to be non-formally smooth. We will use it in the next section to produce concrete examples of this phenomenon. 

\th{necessaria}

If $\Het{i+1}(X,\Zp)$ has non-trivial torsion, then $\ffg{i}{X}$ is not formally smooth.

\eth

\prf

Suppose that there is a non-zero $a$ in the torsion subgroup of $\Het{i+1}(X,\Zp)=\Hcris{i+1}(X)^{F=1}$. We may take $a$ not divisible by $p$. Take $b\in\Acr(C)^{F=p}$ not divisible by $p$. Let $\theta\in\Nyg{}\Hcris{i+1}(X_C)$ the image of $a\otimes b\in\Hcris{i+1}(X)\otimes\Nyg{}\Acr(C)$ via $i_{C,X}$. By design, $(F/p-\iota)(\theta)=0$.

We claim that $a$ maps to a non-zero element $\overline{a}$ in $H^{i+1}(X,\calO_X)$: indeed, $a$ maps to some nonzero element of $\HdR{i+1}(X/k)^{F=1}$, and $F$ acts as zero on the first piece of the Hodge filtration. Therefore $\theta$ maps to the nonzero element $\overline{a}\otimes b$ of 
\begin{equation*}
    \coker\Bigg(i_{X,C}:\Nyg{}\Hcris{i}(X)\otimes_{\W(k)}\Acr(C)\to\Nyg{}\Hcris{i}(X_C)\Bigg)
\end{equation*}
But if $\theta$ were in the image of $\id_X\times f_C$, so would $\overline{a}\otimes b$, which contradicts point $(1)$ of \Cref{L:lemnygaardprodotto}. Therefore $\id_X\times f_C$ is not surjective on $\ker\left(F/p-\iota\right)_{\mathrm{tors}}$, and $\ffg{i}{X}$ is not formally smooth.
\epr

\ssec{esempi}{Varieties with non-formally smooth formal groups}

In this section we take $p=2$ and construct, for any $d\ge2$, a variety $X$ satisfying the following conditions:
\begin{enumerate}
    \item $\fg{i}{X}$ is representable and formally smooth for $i<d$,
    \item $\fg{d}{X}$ is representable but not formally smooth
\end{enumerate}
The first condition ensures that $\fg{d}{X}$ is automatically representable, by \Cref{P:rappartinmazur}. The construction of $X$ is takes inspiration from Igusa's construction \cite{igusa} of a smooth surface with non-reduced Picard variety, which we briefly recall. 

\vspace{.5cm}
\noindent\textbf{Igusa's construction.} Let $E$ be an elliptic curve over $k$ with a nontrivial $2$-torsion point $a$. The automorphism $\sigma$ of $E\times E$ which maps $(x,y)$ to $(x+a,-y)$ is free, so the quotient $X=E\times E/\sigma$ is a smooth surface. To see that the Picard variety of $X$ is non-reduced, we must show that the dimension $g$ of $\Pic_{X/k}$ is strictly smaller than $h^{0,1}=\dim_k(X,\calO_X)$, which is the dimension of its tangent space at $0$. It is easy to show that the map $X\to E/\langle a\rangle$, induced by the first projection, identifies $E/\langle a\rangle$ with the Albanese variety of $X$. Therefore $g=1$. On the other hand, the Hochschild-Serre spectral sequence produces the exact sequence
\begin{equation*}
    0\to k\to H^1(X,\calO_{X})\to H^1(E\times E,\calO_{E\times E})^{\sigma}\to k,
\end{equation*}
where the first copy of $k$ is $H^1(\Z/2,H^0(E\times E,\calO_{E\times E}))$ and the second is $H^0(\Z/2,H^1(E\times E,\calO_{E\times E}))$.  But $\sigma$ acts trivially on $H^1(E\times E,\calO_{E\times E})$, so $h^{0,1}\ge2$. Igusa gives a different, geometric argument for this inequality. 

The first part of the argument can also be replaced by a Hochschild-Serre spectral sequence argument: let $H$ denote either rational $\ell$-adic cohomology for $\ell\ne2$, or rational crystalline cohomology. Recall that if $M$ is a uniquely divisible abelian group with a $\Z/2$-action, then $H^i_{\mathrm{grp}}(\Z/2,M)=0$ for $i>0$, because multiplication by $2$ acts both as zero and as an isomorphism on this group. Therefore, by considering the $E_2$ page of the Hochschild-Serre spectral sequence converging to $H^{p+q}(X)$, one finds that
\begin{equation*}
    H^1(X)=H^1(E\times E)^{\sigma}.
\end{equation*}
Since $g=1/2\dim H^{1}(X)$, a simple computation shows that $g=1$.

\vspace{.5cm}
\noindent\textbf{A computation with Hochschild-Serre.} Let us come back to our problem of constructing varieties $X$ for which $\fg{d}{X}$ is representable, but not formally smooth. Similar to Igusa's variety, we will take $X$ to be a quotient of the form $(E\times Y)/\sigma$ where:

$\bullet$ $Y$ is some variety with an involution $\tau$,
    
$\bullet$ $\sigma$ is the automorphism of $E\times Y$ which maps $(x,y)$ to $(x+a,\tau(y))$,

\noindent plus some conditions on $(Y,\tau)$ given below in \Cref{P:criterioHS}. Note that such an $X$ is always smooth of dimension $\dim Y+1$. To check that $\fg{d}{X}$ is representable and non-formally smooth, we use a Hochschild-Serre spectral sequence argument inspired by the one outlined above. 

\prop{criterioHS}

Suppose that there is an integer $d\ge2$ such that

$(H1)$ $H^i(Y,\calO_Y)=0$ for $0<i<d$,

$(H2)$ $\Hcris{2}(Y)$ is torsion-free,

$(H3)$ there is some nonzero $a\in\Het{d}(Y,\Z_2)$ such that $\tau(a)=-a$.

\noindent Then the following hold

$(C1)$ $H^i(X,\calO_X)=0$ for all $1<i<d$,

$(C2)$ $\Hcris{2}(X)$ is torsion-free,

$(C3)$ $\Het{d+1}(X,\Z_2)_{\mathrm{tors}}\ne0$.

\eprop

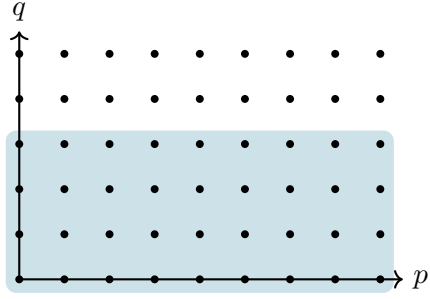
\begin{figure}[h]
    \centering
    \begin{tikzpicture}[scale=0.60]
    \draw[draw=none,fill=SNSBlue,fill opacity=0.2, rounded corners] (-.3,-.3) rectangle (8.3,3.3);
    \foreach \x in {0,...,8}
        \foreach \y in {0,...,5}
            \fill (\x,\y) circle[radius=2.5pt];
    \draw[->,thick] (0,0)--(8.5,0) node[right]{$p$};
    \draw[->,thick] (0,0)--(0,5.5) node[above]{$q$};
    \end{tikzpicture}
    \caption{The region where the spectral sequences of $X$ and $\tilde{E}$ coincide}
    \label{M1}
\end{figure}

\prf
If $H^i(-)$ denotes either crystalline cohomology, étale $2$-adic cohomology, or coherent cohomology, the Hochschild spectral sequence reads
\begin{equation*}
    E^{p,q}_2=H^p_{\mathrm{grp}}(\Z/2,H^q(E\times Y))\implies H^{p+q}(X).
\end{equation*}
First we prove that $(H1)$ implies $(C1)$. Let $\tilde{E}$ denote the quotient $E/\langle a\rangle$. Hypothesis $(H1)$ guarantees that the $E_2$ page for coherent cohomology, with its differentials, coincides with the $E_2$ page of the Hochschild-Serre spectral sequence 
\begin{equation}\label{hsssecoh}
    E_{2}^{p,q}=H^p_{\mathrm{grp}}(\Z/2,H^q(E,\calO_E))\implies H^{p+q}(\tilde{E},\calO_{\tilde{E}})
\end{equation}
in the range $q<d$. Indeed, for $q<d$ the Künneth formula gives $H^q(E\times Y)\simeq H^q(E)$, so the claim is true by functoriality. In \Cref{M1} we draw the picture for $d=4$: in the blue region, the $E_2$ page is the same as for the quotient elliptic curve $\tilde{E}$. 

The $E_2$ page of the Hochschild-Serre spectral sequence \eqref{hsssecoh} is shown in \Cref{paginae2coh}. Note that the differentials drawn in \Cref{paginae2coh} are the only differentials of the whole spectral sequence which can be non-zero. But $\tilde{E}$ is an elliptic curve, so $H^{i}(\tilde{E},\calO_{\tilde{E}})=0$ for $i\ge2$. Therefore all these differentials must be isomorphism, and the $E_3$ page has no non-zero terms outside of $(0,0)$ and $(1,0)$, see \Cref{paginae3coh}.
\begin{figure}[h]
\centering
\begin{minipage}[c]{0.4\linewidth}
\begin{center}
\begin{tikzcd}[column sep = tiny, row sep = tiny]
	{} & 0 & 0 & 0 & 0 \\
	{} & k & k & k & k \\
	{} & k & k & k & k & {} 
	\arrow[from=2-2, to=3-4]
	\arrow[from=2-3, to=3-5]
	\arrow[shorten >=13pt, no head, from=2-4, to=3-6]
\end{tikzcd}
\end{center}
    \caption{the $E_2$ page}
    \label{paginae2coh} 
\end{minipage}
\begin{minipage}[c]{0.4\linewidth}
\begin{center}
\begin{tikzcd}[column sep = tiny, row sep = tiny]
	0 & 0 & 0 & 0 \\
	0 & 0 & 0 & 0 \\
	k & k & 0 & 0
\end{tikzcd}
\end{center}
    \caption{the $E_3$ page}
    \label{paginae3coh}
\end{minipage}
\end{figure}

Now we can prove $(C1)$. By our previous arguments, on page $E_3$ of the spectral sequence 
\begin{equation*}
E_{2}^{p,q}=H^p_{\mathrm{grp}}(\Z/2,H^q(E\times Y,\calO_{E\times Y}))\implies H^{p+q}(X,\calO_{X}),
\end{equation*} 
the only non-zero groups with $q<d$ are $k$ in positions $(0,0)$ and $(1,0)$. Therefore in this region the $E_{\infty}$ page is the same as the $E_3$ page, which tells us that $H^i(X,\calO_X)=0$ for $1<i<d$.

Let's prove $(C2)$. Once again, by the Künneth formula, the two bottom rows of the $E_2$ page of the Hochschild-Serre spectral sequence
\begin{equation}\label{hsssxcris}
E_{2}^{p,q}=H^p_{\mathrm{grp}}(\Z/2,\Hcris{q}(E\times Y)\implies \Hcris{p+q}(X),
\end{equation} 
and the two bottom rows of the $E_2$ page of
\begin{equation}\label{hsssecris}
E_{2}^{p,q}=H^p_{\mathrm{grp}}(\Z/2,\Hcris{q}(E)\implies \Hcris{p+q}(\tilde{E}),
\end{equation}  
are isomorphic. The $E_2$ page of \eqref{hsssecris} is given in \Cref{paginae2}. Once again, the differentials drawn in \Cref{paginae2} are the only differentials of the whole spectral sequence which may be non-trivial. But $\tilde{E}$ is an elliptic curve, so it has torsion-free crystalline cohomology. Therefore the $E_3$ page of \eqref{hsssecris} spectral sequences must be as in \Cref{paginae3}.
\begin{figure}[h]
\centering
\begin{minipage}[c]{0.4\linewidth}
\begin{center}
\begin{tikzcd}[column sep = tiny, row sep = tiny]
	{} & {\W(k)} & 0 & k & 0 & k \\
	{} & {\W(k)^{\oplus2}} & 0 & {k^{\oplus2}} & 0 & {k^{\oplus2}} \\
	{} & {\W(k)} & 0 & k & 0 & k
	\arrow[from=1-2, to=2-4]
	\arrow[from=1-4, to=2-6]
	\arrow[from=2-2, to=3-4]
	\arrow[from=2-4, to=3-6]
\end{tikzcd}
\end{center}
    \caption{the $E_2$ page}
    \label{paginae2} 
\end{minipage}
\begin{minipage}[c]{0.4\linewidth}
\begin{center}
\begin{tikzcd}[column sep = tiny, row sep = tiny]
	2\W(k) & 0 & 0 & 0 \\
	\W(k)\oplus2\W(k) & 0 & 0 & 0 \\
	\W(k) & 0 & 0 & 0
\end{tikzcd}
\end{center}
    \caption{the $E_3$ page}
    \label{paginae3}
\end{minipage}
\end{figure}

By the argument above, the $E_3$ page of \eqref{hsssxcris} is equal to 
\begin{figure}[H]
    \centering
    \begin{tikzcd}[column sep = tiny, row sep = tiny]
	{2\W(k)\oplus M} & \cdots & \cdots & \cdots \\
	{\W(k)\oplus 2\W(k)} & 0 & 0 & 0 \\
	{\W(k)} & 0 & 0 & 0
    \end{tikzcd}
\end{figure}
\noindent in the range $q<3$, where dots stand for some unidentified groups and $M$ is the kernel of the differential in page $E_2$ of \eqref{hsssxcris} mapping $\Hcris{2}(Y)$ to $k^2$. Therefore in the region $p+q<3$ these groups are the same in the $E_{\infty}$ page. Hypothesis $(H2)$ now implies that $\Hcris{2}(X)$ is torsionfree.

Finally, we prove $(C3)$. We claim that $\Het{i}(Y,\Z_2)$ is zero for $0<i<d$. If $a$ is a nonzero element of $\Het{i}(Y,\Z_2)$ which is not divisible by $2$, consider it as an element of $\Hcris{i}(Y)^{F=1}$. It maps to a nonzero element $\overline{a}\in\HdR{i}(Y/k)^{F=1}$, which in turn maps to a nonzero element in $H^i(X,\calO_X)$, because $F$ is zero on the first piece of the Hodge filtration. This contradicts $(H1)$, so indeed $\Het{i}(Y,\Z_2)=0$.

\begin{figure}[h]
\centering
\begin{tikzcd}[column sep = tiny, row sep = tiny]
	0 & 0 & 0 & 0 \\
	{2\Z_2} & 0 & 0 & 0 \\
	{\Z_2} & 0 & 0 & 0
\end{tikzcd}
    \caption{the $E3$ page}
    \label{paginae3et}
\end{figure}
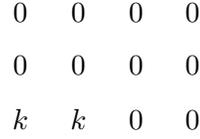

Now we can repeat the argument of the proof of $(C2)$. Namely, in the region $q<d$, the spectral sequence  
\begin{equation}\label{hsssxet}
E_{2}^{p,q}=H^p_{\mathrm{grp}}(\Z/2,\Het{q}(E\times Y,\Z_2)\implies \Het{p+q}(X,\Z_2),
\end{equation}
is the same as the spectral sequence
\begin{equation}\label{hssseet}
E_{2}^{p,q}=H^p_{\mathrm{grp}}(\Z/2,\Hcris{q}(E)\implies \Hcris{p+q}(\tilde{E}),
\end{equation}  
In \Cref{M1} this is depicted as the blue region, for $d=4$. It is easy to see that, similarly to crystalline cohomology, the $E_3$ page of the spectral sequence \cref{hssseet} is as in \Cref{paginae3et}, with zero differentials.
\noindent Therefore, on the $E_{\infty}$ page of \cref{hsssxet} we have $E_{\infty}^{i,d+1-i}=0$ for $i\ge2$ and $E_{\infty}^{1,d}=H^1(\Z/2,\Het{d}(Y,\Z_2))$, which has nonzero torsion as a consequence of $(H3)$. Thus $\Het{d+1}(X,\Z_2)$ contains $H^1(\Z/2,\Het{d}(Y,\Z_2))$ and also has nonzero torsion. This proves $(C3)$.
\epr

\cor{criterioesempio}

If such a $Y$ exists, the functors $\fg{i}{X}$ are all representable and formally smooth for $i<d$, and $\fg{d}{X}$ is representable but not formally smooth.

\ecor{}

\prf
Since $\Hcris{2}(X)$ is torsionfree, $\fg{1}{X}$ is formally smooth by \Cref{T:sufficiente}. Recall (\Cref{P:compconn}) that the tangent space to $\fg{i}{X}$ is a quotient $H^i(X,\calO_X)$, therefore $\fg{i}{X}=0$ for $1<i<d$. Artin--Mazur's representability criterion \Cref{P:rappartinmazur} then implies that $\fg{d}{X}$ is representable, and finally thanks to \Cref{T:necessaria} we have that $\fg{d}{X}$ is not formally smooth. 
\epr

\vspace{.25cm}
\noindent\textbf{A constructive example for $d=2$.} We claim that if $Y$ is an ordinary K3 surface with an Enriques involution $\tau$, then the hypothesis of \Cref{P:criterioHS} are satisfied. The crystalline cohomology of a K3 surface $Y$ is always torsionfree and $\Hcris{1}(S)=0$ - see \cite[II.7.2]{illusiedrw} - so the first condition is met. If $Y$ is also ordinary we have a decomposition of crystals
\begin{equation}
    \Hcris{2}(Y)\simeq\W(0)\oplus\W(1)^{20}\oplus\W(2),
\end{equation}
where $\W(i)$ is the rank one crystal $x\cdot\W(k)$ with $F(x)=p^ix$. Therefore $\Het{2}(Y,\Z_2)\simeq\Z_2$. Suppose finally that $Y$ has a fixed-point-free involution $\tau$. Then $\tau$ acts as $\pm\id$ on $\Het{2}(Y,\Z_2)$. The quotient $\tilde{Y}=Y/\tau$ is an Enriques surface, for which it is known \cite[II.7.3]{illusiedrw} that $\Het{2}(\tilde{Y},\Q_2)=0$. Since $\Het{2}(\tilde{Y},\Q_2)=\Het{2}(Y,\Q_2)^{\tau}$, the action must be by $-\id$.

It remains to find such a $Y$, which does not seem to appear in existing literature. The following example was constructed with help from Yuya Matsumoto: let $A=E'\times E'$ where $E'$ is an ordinary elliptic curve, and let $b$ be its nontrivial $2$-torsion point. We take $Y$ to be the Kummer surface associated to $A$, which is a $K3$ surface of Picard rank $20$, see \cite{shioda} for details. By design, $Y$ admits an elliptic fibration, therefore by a theorem of Artin \cite[Theorem 1.7]{artinssk3} it is not supersingular, so it must be ordinary (see \cite[II.7.2]{illusiedrw}). Consider the involution $\tau'$ of $A$, which maps $(x,y)$ to $(x-b,-y+b)$. If $(x,y)\in A\times A$, the equalities $(x,y)=(x-b,-y+b)$ and $(-x,-y)=(x-b,-y+b)$ are both impossible. Therefore $\tau'$ induces a fixed-point free involution $\tau$ of $Y$, which is what we wanted to show.

\prop{esempietto}
    Let $X=E\times Y/\sigma$ with $Y$ as in the previous paragraph. Then $\Pic_{X/k}$ is smooth, and $\widehat{Br}(X)$ is representable but not formally smooth.\qed
\eprop

The variety $X$ we just defined is a family of ordinary Enriques surfaces fibred over the elliptic curve $E/\langle a\rangle$.

\vspace{.5cm}
\noindent\textbf{A nonconstructive example for any $d\ge2$.} The construction of this example is inspired by an argument of Koblitz \cite{koblitz} on the Hasse-Witt matrix of complete intersections. Let $n\ge3$, and let $Y\subseteq\mathbb{P}^n_k$ be a hypersurface of degree $d_1$, where $d_1\ge1$. There is a short exact sequence
\begin{equation*}
    0\to\calO_{\mathbb{P}^n_k}(-d_1)\to\calO_{\mathbb{P}^n_k}\to\calO_Y\to0
\end{equation*}
of coherent sheaves on $\mathbb{P}^n_k$, which shows
\begin{equation*}
    H^i(Y,\calO_Y)=0\textnormal{ for }0<i<n-1,\quad H^{n-1}(Y,\calO_Y)\simeq H^{n}(\mathbb{P}^n_k,\calO_{\mathbb{P}^n_k}(-d_1)),
\end{equation*}
which is zero for $d_1<n+1$ and non-zero otherwise. Recall that for a projective $k$-variety $Z$ of dimension $m$, the Hasse-Witt matrix of $Z$ is the semilinear endomorphism $F$ of $H^{m}(Z,\calO_Z)$ induced by the absolute Frobenius of $Z$. We say the Hasse-Witt matrix of $Z$ is invertible if $F$ is bijective. In our setting, we will thus be concerned with the action of $F$ on $H^{n-1}(Y,\calO_Y)$.

Let us recall Koblitz's proof that the general hypersurface of degree $d$ in $\mathbb{P}^n$ has invertible Hasse-Witt matrix. Consider the coordinates $x_0,\dots,x_n$ on $\mathbb{P}^n_k$. Hypersurfaces of degree $d$ are parametrized by homogeneous forms of degree $d$ in the variables $x_0,\dots,x_n$, up to scalar multiplication. The space of such parameters is naturally identified with a projective space $\mathbb{P}^N$, where $N={{d+n}\choose{d}}$, in the sense that there is a flat and proper family $\mathcal{H}\to\mathbb{P}^N$ whose fibre over $[y]$ is the hypersurface defined by the equation $y=0$. Koblitz shows that hypersurfaces with invertible Hasse-Witt matrix cut out an open subvariety of $\mathbb{P}^N$. To exclude that this open set is empty, he shows that the union of $d$ hyperplanes intersecting properly has invertible Hasse-Witt matrix, thus completing the proof.

We will need a variant of this: let $\tau$ be the involution of $\mathbb{P}^N$ exchanging $x_0$ and $x_1$. There is a linear subspace $\mathbb{P}^{m}$ of $\mathbb{P}^N$ which parametrizes hypersurfaces preserved by $\tau$: in coordinates, it consists of the linear forms which are invariant (up to scalar multiplication) under exchanging $x_0$ and $x_1$. Thus, the family of all hypersurfaces over $\mathbb{P}^N$ restricts to a flat family of hypersurfaces $\mathcal{H}'\to\mathbb{P}^m$.

\lem{}

The generic fibre of the family $\mathcal{H}'$ is smooth.

\elem

\prf

It suffices to produce one smooth hypersurface of degree $d_1$ preserved by $\tau$. If $d_1$ is odd we may take the Fermat hypersurface of degree $d_1$ in $\mathbb{P}^n$. If $d_1$ is even, consider
\begin{equation*}
    f(x_0,\dots,x_n)=x_0^{d_1/2}x_1^{d_1/2}+\varepsilon+x_2^{d_1}+\sum_{i=2}^{n-1}x_ix_{i+1}^{d_1-1}+x_n(x_0^{d_1-1}+x_1^{d_1-1}),
\end{equation*}
where $\varepsilon$ is defined as follows: 
\begin{equation*}
    \varepsilon=\begin{cases}
        0, \quad &d_1/2\textnormal{ is odd}\\
        x_2^{d_1-2m}(x_0x_1)^{m}, \quad &d_1=2^sm,\textnormal{ $m$ odd, $s\ge2$}
    \end{cases}
\end{equation*}
The equation $f=0$ defines a hypersurface of degree $d_1$ in $\mathbb{P}^n$ preserved by $\tau$. We claim it is non-singular: let $a=(a_0,\dots,a_n)\in k^{n+1}$ satisfy
\[
f(a)=\frac{\partial f}{\partial x_0}(a)=\dots=\frac{\partial f}{\partial x_n}(a)=0.
\]
The equations $\frac{\partial f}{\partial x_i}(a)=0$ for $2\le i\le n-1$ give inductively $a_j=0$ for $j\ge3$. Then, $\frac{\partial f}{\partial x_n}=0$ yields 
\begin{equation}\label{eqb}
    a_0^{n-1}=a_1^{n-1}.
\end{equation} 
Now, we must distinguish two cases. First suppose that $d_1/2$ is odd. Then equation $\frac{\partial f}{\partial x_0}(a)=0$ implies that either $a_1=0$ or $a_0=0$, which combined with \eqref{eqb} gives $a_1=a_0=0$. Then $f(a)=0$ finally implies $a_2=0$, i.e. $a_i=0$ for all $i$.

If $d_1/2$ is even, equation $\frac{\partial f}{\partial x_0}(a)=0$ implies that either $a_0=0$ or $a_1=0$ or $a_2=0$. In the first two cases, we have once again $a_0=a_1=0$ and thus $a_2=0$ as above. In the third case, $f(a_0,\dots,a_n)=0$ implies $a_0=0$ or $a_1=0$, and \eqref{eqb} shows $a_0=a_1=0$. In either case $a_i=0$ for all $i$. 

This shows that $\{f=0\}$ is nonsingular, and the lemma is proved.
\epr

\lem{esistenzaiper}

There exists a smooth hypersurface $Y\subseteq\mathbb{P}^n$ of degree $d_1$ with invertible Hasse-Witt matrix.

\elem

\prf
By the proof of \cite[I, Lemma 4]{koblitz}, if the Hasse-Witt matrix is invertible for one fibre of $\mathcal{H}'$, then the same holds for a general fibre. Therefore to prove the lemma it will suffice to exhibit a possibly singular hypersurface, preserved by $\tau$, with invertible Hasse-Witt matrix. By \cite[I, Theorem 3]{koblitz}, a union of $d_1$ hyperplanes in general position has invertible Hasse-Witt matrix. To conclude, take such a collection of hyperplanes preserved by $\tau$.
\epr

For the next result, we take $n=d+1$ and $d_1=d+2$, thus $Y$ is a Calabi-Yau variety of dimension $d$.

\th{esempio} Let $Y$ be a hypersurface as in \Cref{L:esistenzaiper}. Then, if $X=E\times Y/\sigma$, the functors $\fg{i}{X}$ are all representable and formally smooth for $i<d$, and $\fg{d}{X}$ is representable but not formally smooth. \qed

\eth

\prf
Any smooth hypersurface satisfies the hypotheses $(H1)$ and $(H2)$ of \Cref{P:criterioHS} by the Lefschetz theorem. We claim that $(H3)$ is also verified for $(Y,\tau)$ - then the conclusion will follow from \Cref{C:criterioesempio}.

Consider $\mathbb{P}^{d+1}_{\W(k)}$ with coordinates $y_0,\dots,y_{d+1}$ lifting $x_0,\dots,x_{d+1}$. We claim that we can lift $Y\subseteq\mathbb{P}^{d+1}$ to $\tilde{Y}\subseteq\mathbb{P}^{d+1}_{\W(k)}$, in such a way that $\tilde{Y}$ is preserved by $\tilde{\tau}$, the automorphism of $\mathbb{P}^{d+1}_{\W(k)}$ exchanging $y_0$ and $y_1$. Indeed, consider a non-singular hypersurface of degree $s$,
\begin{equation*}
S=\Bigg\{\sum_{|I|=s} a_Ix^I=0\Bigg\}\subseteq\mathbb{P}^n_k
\end{equation*}
Then we may define the subscheme
\begin{equation*}
\tilde{S}=\Bigg\{\sum_{|I|=s} [a_I]y^I=0\Bigg\}\subseteq\mathbb{P}^n_{\W(k)},
\end{equation*}
where $[a]$ denotes the Teichmüller lift of $a\in k$. Both the special fibre and the generic fibre of $\tilde{S}$ are smooth hypersurfaces of the same dimension, therefore $\tilde{S}$ is flat over $\Spec(\W(k))$. Then, if we define $\tilde{Y}$ by this procedure, the hypersurface $\tilde{Y}$ is fixed by $\tilde{\tau}$. Indeed, being fixed by $\tau$ is equivalent to certain coefficients of the defining equations being equal, and this condition is preserved upon taking Teichmüller lifts. 
\begin{equation*}
    \Hcris{d}(Y)\simeq\HdR{d}(\tilde{Y}/\W(k))
\end{equation*}
and furthermore a surjection
\begin{equation*}
    \HdR{d}(\tilde{Y}/\W(k))\onto{}H^d(\tilde{Y},\calO_{\tilde{Y}})\simeq H^{d+1}(\mathbb{P}^{d+1}_{\W(k)},\omega_{\mathbb{P}^{d+1}_{\W(k)}}).
\end{equation*}
The rightmost group is a free $\W(k)$-module of rank one: using \v{C}ech cohomology with respect to the standard covering of $\mathbb{P}^{d+1}_{\W(k)}$, a generator of this group is the cocycle
\begin{equation*}
    \alpha=\frac{y_{d+1}^{d+1}}{y_0\dots y_{d}}d\left(\frac{y_0}{y_{d+1}}\right)\wedge\dots\wedge d\left(\frac{y_{d}}{y_{d+1}}\right),
\end{equation*}
see \cite[III, Remark 7.1.1]{hartshorne}, for example. We see that $\tilde{\tau}(\alpha)=-\alpha$. Therefore $\tilde{\tau}$ acts as $-\id$ on $H^d(\tilde{Y},\calO_{\tilde{Y}})$.

Let $M_0$ be the unit-root sub-$F$-crystal of $\Hcris{d}(Y)$. The surjection $\Hcris{d}(Y)\to H^{d}(Y,\calO_Y)$ induces a surjection $M_0\to\HdR{d}(Y/k)^{ss}$, where the latter is the semisimple part of $\HdR{d}(X/k)$, i.e. the largest subspace of $\HdR{d}(Y/k)$ on which $F$ acts bijectively. Now $\HdR{d}(Y/k)\to H^{d}(X,\calO_X)$ is surjective and compatible with $F$, so it is surjective upon taking semisimple parts. Therefore, we obtain a surjection $M_0\onto{} H^{d}(Y,\calO_Y)^{ss}$, inducing a surjection $M_0^{F=1}\onto{} H^{d}(Y,\calO_Y)^{F=1}$, so there is $a\in M_0^{F=1}=\Hcris{d}(Y)^{F=1}$ which maps to a non-zero element of $H^{d}(Y,\calO_Y)$. Then $a$ also maps to a non-zero element of $H^{d}(\tilde{Y},\calO_{\tilde{Y}})$. By the above paragraph $\tau(a)=-a$ and $(H3)$ holds, which is what we wanted to show.
\epr

\ssec{pruffafinale}{Proof of \Cref{T:equivalenteintro}}

For this final section, $p$ is any prime number, and $X$ is a smooth and proper variety over $k$. We start by recalling the essential properties of Witt vector cohomology following the paper \cite{serremexico} where it was introduced.

Let $n\ge1$. We can regard the truncated Witt vectors $\mathcal{W}_n$ as an affine group scheme over $k$, by defining its functor of points to be $\Spec(A)\mapsto\W_n(A)$. Each $\mathcal{W}_n$ is a successive extension of copies of $\Ga$, therefore its cohomology on $X$ can be computed indifferently in the Zariski, étale or fppf topology. It also follows by dévissage that $H^i(X,\mathcal{W}_n)$ is a finite length $\W(k)$-module for any $i\ge0$. Thus we can define the $q$-th Witt vector cohomology group either as a projective limit along the restriction maps $\mathcal{W}_{n+1}\to\mathcal{W}_n$,
\begin{equation*}
    H^i(X,\mathcal{W})=\varprojlim_n H^i(X,\mathcal{W}_n),
\end{equation*}
or as the cohomology of the complex $\Rlim\Rflat(X,\mathcal{W}_n)$.

The system of the $\mathcal{W}_n$ is endowed with morphisms of algebraic groups $F$ and $V$, respectively the Frobenius and the Verschiebung, which satisfy the usual relations. The short exact sequences
\[
0\to \mathcal{W}_{n-1}\xto{}{V}\mathcal{W}_n\to\Ga\to0
\]
yield a long exact sequence in cohomology, and upon taking the projective limit we get a long exact sequence
\begin{equation}\label{lesversch}
\dots\to H^i(X,\mathcal{W})\xto{}{V} H^i(X,\mathcal{W})\to H^{i}(X,\calO_X)\to H^{i+1}(X,\mathcal{W})\to\dots
\end{equation}

This is what we need of the classical theory to state and prove this section's result. In order to relate this to our previous discussion, we give a description of these groups with quasisyntomic descent, via a ``modified Nygaard filtration'' construction. The first remark we make is that $\Rflat(-,\mathcal{W}_n)$ satisfies pro-fppf descent, and in particular quasisyntomic descent, as in \cite[Remark 7.2.4]{bhattlurie}.

\defe{}
    For a quasisyntomic $k$-scheme $Y$, there is a natural map
    \begin{equation*}
        \Rcris(Y/\Zp)\to\Rflat(Y,\mathcal{W})
    \end{equation*}
    which is defined as follows: if $Y=\Spec(A)$ for some qrsp $k$-algebra $A$, it corresponds to the surjection $\Acr(A)\to W(A)$, coming from $\W(\perf{A})\to\W(A)$. Then we extend this definition to all quasisyntomic $k$-algebras via descent. We thus have for all $n$ a natural map
    \begin{equation}\label{defefunoen}
        \Rcris(Y/\Zp)\to\Rflat(Y,\mathcal{W}_n).
    \end{equation}
    
    $(1)$ Denote by $\Nygm{n}\Rcris(Y)$ the fibre of \eqref{defefunoen}, and by $\Nygm{n}\Hcris{i}(Y)$ its $i$-th cohomology group. 
    
    $(2)$ If $Y=\Spec(A)$ for some qrsp $k$-algebra $A$, the complex $\Nygm{n}\Rcris(Y)$ is identified with the kernel of $\Acr(A)\to\W_n(A)$ concentrated in degree $0$, which we call $\Nygm{n}\Acr(A)$.
\edefe

An interesting feature of this complex is that it can be described, for eqrsp $k$-algebras, purely in terms of the divisibility of Frobenius. We make this precise in the following lemma - we believe it should hold for any qrsp $k$-algebra, but we do not need it here.
 
\lem{descrizionenygmod}
Let $A$ be eqrsp, as in \eqref{eqrsp}. We have
\begin{equation}\label{descrizionenygmod}
    \Nygm{n}\Acr(A)=\{a\in\Acr(A)\textnormal{ s.t. }p^i|F^i(a)\textnormal{ for all } 1\le i\le n\}.
\end{equation}
\elem

\prf
In the notation of \Cref{L:descrizioneacris}, the right-hand group is the ideal of $\Acr(A)$ topologically generated by elements of the form $x^{<\alpha>}, \alpha\ge1$, and of the form $p^{m}x_i^{<p^{-m}>}, m\ge0$. According to \cite[2.5.3]{drinfeldacris}, the kernel of $\W(\perf{C})\to\W(C)$ is topologically generated by elements of the form $p^mx_i^{p^{-m}}$, so equality \eqref{descrizionenygmod} holds.
\epr

Therefore, if $Y$ is the spectrum of an eqrsp algebra, there are natural maps 
\begin{equation*}
    (F^i/p^i):\Nygm{m}\Rcris(Y)\to \Nygm{m-i}\Rcris(Y)
\end{equation*}
for any $1\le i\le m$, where we write $\Nygm{0}\Rcris(Y)$ for $\Rcris(Y)$. We can extend the definition of these maps to smooth schemes by descent.

\lem{commutation}
If $Y$ is a smooth $k$-scheme, and $n\ge 1$, consider the commutative diagram
\begin{equation*}
    \begin{tikzcd}
	{\Nygm{n+1}\Rcris(Y)} & {\Nyg{}\Rcris(Y)} \\
	{\Nygm{n}\Rcris(Y)} & {\Rcris(Y)}
	\arrow[from=1-1, to=1-2]
	\arrow["{F/p}"', from=1-1, to=2-1]
	\arrow["{F/p}"', from=1-2, to=2-2]
	\arrow[from=2-1, to=2-2]
\end{tikzcd}
\end{equation*}
in $\D(\Z_p)$. It induces a quasi-isomorphism
\begin{equation*}
    \Cone\left(\Nygm{n+1}\Rcris(Y)\to\Nyg{}\Rcris(Y)\right)\simeq\Rflat(Y,\mathcal{W}_n).
\end{equation*}
\elem

\prf
By descent, we reduce to the following statement: let $A$ be an eqrsp $k$-algebra. If $M$ denotes the quotient of $\Nyg{}\Acr(A)$ by its subgroup $\Nygm{n+1}\Acr(A)$, then $F/p:\Nyg{}\Acr(A)\to\Acr(A)$ induces an isomorphism
\begin{equation*}
    \Nyg{}\Acr(A)/\Nygm{n+1}\Acr(A)\simeq\Acr(A)/\Nygm{n}\Acr(A)\simeq\W_n(A).
\end{equation*} 
This is straightforward to check, using the explicit descriptions of these groups given in \Cref{L:descrizionenygmod}.
\epr

\lem{successionefsup}
Let $a\in\Nyg{}\Hcris{i}(X)$ and $n\ge1$. The following are equivalent.

$(1)$ the element $a$ lifts to $\Nygm{n}\Hcris{i}(X)$.

$(2)$ there is a sequence $a_1,\dots,a_n\in\Nyg{}\Hcris{i}(X)$ with $a_1=a$ and $(F/p)(a_j)=\iota(a_{j+1})$ for $1\le j\le n-1$.
\elem

\prf
Suppose $(1)$ holds, so $a$ is the image of some $a'\in \Nygm{n}\Hcris{i}(X)$. Condition $(2)$ is satisfied by taking $a_j$ to be the image of $\left(F^{j-1}/p^{j-1}\right)(a')\in \Nygm{n-j}\Hcris{i}(X)$ in $\Nyg{}\Hcris{i}(X)$. 

Now suppose $(2)$ holds. We prove by induction that $a_{n-j}$ lifts to $\Nygm{j+1}\Hcris{i}(X)$. For $j=0$ this is clear. The induction step is proved by taking cohomology in the diagram of \Cref{L:commutation}.
\epr

\th{equivalente}
The following are equivalent.

$(1)$ the formal group $\fg{i}{X}$ is formally smooth.

$(2)$ the map $H^{i}(X,\mathcal{W})\to H^{i}(X,\calO_X)$ is surjective.

$(3)$ the Verschiebung $V$ acting on $H^{i+1}(X,\mathcal{W})$ is injective.
\eth

\prf[Proof of \Cref{T:equivalente}]

We start with some preliminary remarks. First, it is clear from the sequence \eqref{lesversch} that conditions $(2)$ and $(3)$ are equivalent. Moreover that same sequence was obtained as the projective limit of sequences involving the cohomology of $\mathcal{W}_n$, from which we see that condition $(2)$ is equivalent to the natural maps $H^{i}(X,\mathcal{W}_n)\to H^{i}(X,\calO_X)$ being surjective for all $n\ge 1$.

Using the definition of $\Nygm{n}\Rcris(X)$, we have for each $n$ a commutative diagram
\begin{equation*}
\begin{tikzcd}
	{\Hcris{i}(X)} & {H^{i}(X,\calO_X)} & {\Nyg{}\Hcris{i+1}(X)} & {\Hcris{i+1}(X)} \\
	{\Hcris{i}(X)} & {H^i(X,\mathcal{W}_n)} & {\Nygm{n}\Hcris{i+1}(X)} & {\Hcris{i+1}(X)}
	\arrow[from=1-1, to=1-2]
	\arrow[from=1-2, to=1-3]
	\arrow[from=1-3, to=1-4]
	\arrow[shift right, no head, from=1-4, to=2-4]
	\arrow[shift right, no head, from=2-1, to=1-1]
	\arrow[shift left, no head, from=2-1, to=1-1]
	\arrow[from=2-1, to=2-2]
	\arrow[from=2-2, to=1-2]
	\arrow[from=2-2, to=2-3]
	\arrow[from=2-3, to=1-3]
	\arrow[from=2-3, to=2-4]
	\arrow[shift right, no head, from=2-4, to=1-4]
\end{tikzcd}
\end{equation*}
from which we see that $H^i(X,\mathcal{W}_n)\to H^i(X,\calO_X)$ is surjective if and only if $V^i_X\subseteq\Nyg{}\Hcris{i+1}(X)$ is in the image of $\Nygm{n}\Hcris{i+1}(X)\to\Nyg{}\Hcris{i+1}(X)$. Hence, by \Cref{L:successionefsup}, we find that condition $(2)$ is equivalent to the following:

$(2')$ if $a\in V^i_X\subseteq\Nyg{}\Hcris{i+1}(X)$, for any $n\ge0$ there is a sequence $a_1,\dots,a_n\in\Nyg{}\Hcris{i+1}(X)$ with $a_1=a$ and $(F/p)(a_j)=\iota(a_{j+1})$.

Note that $F/p:\Nyg{}\Hcris{i+1}(X)\to\Hcris{i+1}(X)$ is injective modulo torsion, so condition $(2')$ is equivalent to:

$(2'')$ if $a\in V^i_X\subseteq\Nyg{}\Hcris{i+1}(X)$, for any $n\ge0$ there is a sequence $a_1,\dots,a_n\in\Nyg{}\Hcris{i+1}(X)_{\mathrm{tors}}$ with $a_1=a$ and $(F/p)(a_j)=\iota(a_{j+1})$.

Summing up, it remains to show that condition $(1)$ is equivalent to condition $(2'')$. By \Cref{C:criteriofinale}, the former is equivalent to $\id_X\times f_C$ acting surjectively on $\ker\left(F/p-1\right)_{\mathrm{tors}}\subseteq\Nyg{}\Hcris{i+1}(X_C)$, and using \Cref{L:lemnygaardprodotto} and \eqref{formulafsup} we can describe this group explicitly. Namely, we see that an element of $\Nyg{}\Hcris{i+1}(X)$, written as in \eqref{formaelt}, lies in $\ker\left(F/p-1\right)_{\mathrm{tors}}$ if and only if the following conditions are met:
\begin{enumerate}
    \item[(K1)] each $a_{\alpha}$ and each $b_{\alpha}$ is a torsion element of the group it belongs to, and all sums are finite.
    \item[(K2)] if $\alpha<1$ we have $(F/p)(a_{\alpha/p})=\iota(a_{\alpha})$.
    \item[(K3)] if $1\le\alpha<p$ we have $(F/p)(a_{\alpha/p})=b_{\alpha}+\iota(a_{\alpha})$.
    \item[(K4)] if $p\le\alpha$ we have $p^{\lfloor\alpha/p\rfloor}(F/p)(a_{\alpha/p})+p^{\lfloor\alpha/p\rfloor-1}F(b_{\alpha/p})=b_{\alpha}+\iota(a_{\alpha})$.
\end{enumerate}

Suppose $(1)$ holds, so that $\id_X\times f_C$ acts surjectively on the set of such elements. Then any such element must satisfy $b_{\alpha}=0$ for all $\alpha$, by \eqref{formulafC}. If $a\in V^i_X$, we have $a\otimes x^{<1/p>}+(F/p)(a)\otimes x^{<1>}\in\ker\left(F/p-1\right)_{\mathrm{tors}}$. Take $a'$ in the same group such that $f_C^n(a')=a$. In the notation of \eqref{formaelt}, $a'$ has monomials of the form $a_{n-j}\otimes x^{<p^{-j-1}>}$, for $0\le j\le n-1$, which satisfy the relations $(F/p)(a_{j})=\iota(a_{j+1})$ and $a_{1}=a$. Therefore condition $(2'')$ is verified.

Now suppose $(2'')$ holds. By \Cref{L:successionefsup}, for every $m\ge1$ and every $a\in V^{i}_X$, there is an element of $\Nygm{m}\Hcris{i+1}(X)_{\mathrm{tors}}$ which maps to $a$ via the natural map $\Nygm{m}\Hcris{i+1}(X)\to\Nyg{}\Hcris{i+1}(X)$.

Let $e_{i+1}$ be the $p$-exponent of $\Hcris{i+1}(X)$. Let $(\overline{a},\overline{\alpha})\in \Nygm{m}\Hcris{i+1}(X)_{\mathrm{tors}}\times\Z_+[1/p]$ be a pair of elements such that the image of $\overline{a}$ in $\Nyg{}\Hcris{i+1}(X)$ lies in $V^i_X$, and $p^{m-1}\overline{\alpha}\ge e_{i+1}$. We define an element $\theta_{\overline{a},\overline{\alpha}}$ in $\ker\left(F/p-1\right)_{\mathrm{tors}}$ as follows:
\begin{equation*}
    \theta_{a,\alpha}=\sum_{i=0}^{m-1} a_i\otimes x^{<p^i\overline{\alpha}>},\quad\quad a_i=p^{\sum_{j=0}^{i-1}\lfloor p^j\overline{\alpha}\rfloor}(F^i/p^i)(\overline{a}).
\end{equation*}
Note that if $p^m\overline{\alpha}\ge e_{i+1}$, then $\theta_{\overline{a},\overline{\alpha}/p}$ is well-defined, and $(\id_X\times f_C)(\theta_{\overline{a},\overline{\alpha}/p})=\theta_{\overline{a},\overline{\alpha}}$. Therefore, if we show that $\ker\left(F/p-1\right)_{\mathrm{tors}}$ is generated by such elements, we will be done.

Let $y\in\ker\left(F/p-1\right)_{\mathrm{tors}}$. By adding elements of the type $\theta_{\overline{a},\overline{\alpha}}$ to $y$, with $p^m\overline{\alpha}\ge e_{i+1}$, we may arrange that the smallest $\alpha$ such that $a_{\alpha}\ne0$, say $\alpha_0$, is greater or equal than $1$. Indeed, if $\alpha_0<1$, by condition K2 we have $a_{\alpha_0}\in V^i_X$, so by the running hypothesis we may choose an appropriate lift to $\Nygm{m}\Hcris{i+1}(X)$ for $m>>0$. We claim that $b_{\alpha}=0$ for all $\alpha\ge1$. Condition $K3$ shows that $b_{\alpha}+\iota(a_{\alpha})=0$ for $1\le\alpha<p$. But $b_{\alpha}\in\Nyg{}\Hcris{i+1}(X)$ if and only if $b_{\alpha}=0$, so we must have $b_{\alpha=0}$ for $1\le\alpha<p$. We may now use condition K4 to repeat the argument, and show by induction that $b_{\alpha}=0$ for all $p^{i}\le\alpha<p^{i+1}$. Thus $b_{\alpha}=0$ for all $\alpha$.

Now we have $\alpha_0\ge1$. Condition K3 (or K4) shows that $\iota(a_{\alpha_0})=0$. Thus $a_{\alpha_0}=0$, because $a_{\alpha}$ is determined mod $V^i_X$ for $\alpha\ge1$. By our definition of $\alpha_0$ we must have $a_{\alpha}=0$ for all $\alpha$, hence $y=0$, and we are done.

\epr

\vspace{.5em}
\hrule
\vspace{1em}
\eu

\end{document}